\documentclass[12pt]{amsart}
\usepackage{a4wide}
\usepackage{amssymb}

\usepackage[whole,autotilde]{bxcjkjatype}

\usepackage[mathscr]{eucal}

\usepackage{ifmtarg}
\usepackage{version}
\usepackage{tikz}
\usetikzlibrary{arrows,decorations.pathmorphing,backgrounds,positioning,fit,petri}
\usetikzlibrary{decorations.markings}

\usepackage{url}
\usepackage[
pdftex,
bookmarks=true,
colorlinks=true,
citecolor=deepjunglegreen,
filecolor=deepjunglegreen,
debug=true,
naturalnames=true,
unicode,
pdfnewwindow=true]{hyperref}

\usepackage{color}
\definecolor{deepjunglegreen}{rgb}{0.0, 0.29, 0.29}

\newenvironment{NB}{
\color{red}{\bf NB}. \footnotesize
}{}
\newenvironment{NB2}{
\color{blue}{\bf NB2}. \footnotesize
}{}
\excludeversion{NB}
\excludeversion{NB2}

\newenvironment{aenume}{%
  \begin{enumerate}%
  }{\end{enumerate}}
\newcommand{\defeq}{\overset{\scriptstyle\mathrm{def.}}{=}}
\newcommand{\CC}{{\mathbb C}}
\newcommand{\ZZ}{{\mathbb Z}}

\newcommand{\proj}{{\mathbb P}}

\newcommand{\GL}{\operatorname{GL}}

\newcommand{\grpSp}{\operatorname{\rm Sp}}

\newcommand{\algsl}{\operatorname{\mathfrak{sl}}} 

\newcommand{\Spec}{\operatorname{Spec}\nolimits}
\newcommand{\Proj}{\operatorname{Proj}\nolimits}
\newcommand{\End}{\operatorname{End}}
\newcommand{\Hom}{\operatorname{Hom}}

\renewcommand{\MR}[1]{}

\newcommand{\dslash}{/\!\!/}
\newcommand{\vin}[1]{\operatorname{i}(#1)} 
\newcommand{\vout}[1]{\operatorname{o}(#1)} 
\newcommand{\bM}{\mathbf M}
\newcommand{\bN}{\mathbf N}

\newcommand{\bG}{\mathbf G}

\newcommand{\tslash}{/\!\!/\!\!/}
\newcommand{\tslabar}{\mathbin{
\setbox0=\hbox{/\!\!/\!\!/}\rule[0.4\ht0]{\wd0}{.3\dp0}\kern-\wd0\box0}}

\newcommand{\Gr}{\mathrm{Gr}}

\newcommand{\cR}{\mathcal R}

\newcommand{\cT}{\mathcal T}
\newcommand{\cK}{\mathcal K}
\newcommand{\cO}{\mathcal O}

\newcommand{\scP}{\mathscr P}

\makeatletter
\newcommand{\cA}[1][{}]{%
  \@ifmtarg{#1}%
  {\mathcal A}
  {\mathcal A(#1)}
}
\newcommand{\cAh}[1][{}]{%
  \@ifmtarg{#1}%
  {\mathcal A_\hbar}
  {\mathcal A_\hbar(#1)}
}
\makeatother

\newcommand{\fA}{\mathfrak A}

\newcommand{\ft}{\mathfrak t}

\newcommand{\po}{\ar@{}[dr]|{\text{\pigpenfont R}}}
\newcommand{\pb}{\ar@{}[dr]|{\text{\pigpenfont J}}}
\newcommand{\pp}{\ar@{}[dr]|{\text{\pigpenfont P}}}

\newcommand{\cM}{\mathcal M}

\newcommand{\BA}{{\mathbb{A}}}

\newcommand{\BG}{{\mathbb{G}}}

\newcommand{\bT}{{\mathbf{T}}}
\newcommand{\bp}{{\mathbf{p}}}
\newcommand{\bq}{{\mathbf{q}}}

\newcommand{\oW}{\overline{\mathcal{W}}{}}

\newcommand{\fg}{{\mathfrak{g}}}

\newcommand{\fM}{{\mathfrak{M}}}

\newcommand{\fL}{\mathfrak L}
\newcommand{\fT}{\mathfrak T}

\newcommand{\Perv}{\operatorname{Perv}}
\newcommand{\IC}{\operatorname{IC}}

\DeclareSymbolFont{symbolsC}{U}{pxsyc}{m}{n}
\SetSymbolFont{symbolsC}{bold}{U}{pxsyc}{bx}{n}
\DeclareFontSubstitution{U}{pxsyc}{m}{n}
\DeclareMathSymbol{\medcirc}{\mathbin}{symbolsC}{7}

\theoremstyle{definition}

\newtheorem{Theorem}{定理}[subsection]

\newtheorem{Lemma}[Theorem]{補題}
\newtheorem{Proposition}[Theorem]{命題}

\newtheorem{Definition}[Theorem]{定義}

\newtheorem{Conjecture}[Theorem]{予想}
\newtheorem{Remark}[Theorem]{注}

\newcommand{\secref}[1]{\S\ref{#1}}
\newcommand{\subsecref}[1]{\ref{#1}}

\begin{document}



\title[クーロン枝と幾何学的佐武対応]
{超対称性ゲージ理論のクーロン枝の数学的定義と\\
Kac-Moodyリー環の幾何学的佐武対応}
\author{中島　啓}
\maketitle



\section*{はじめに}

筆者とBraverman，Finkelbergは，{\bf クーロン枝 (Coulomb branch)}とよば
れる新しい代数多様体のクラスと，その非可換変形(量子化)を導入し
た \cite{2015arXiv150303676N,2016arXiv160103586B}．
%
%
この研究は，もともと理論物理学の超対称性ゲージ理論において研究されてい
たクーロン枝に，数学的に厳密な定義を与える，という動機に基づいている．
そこで得た定義は，ホモロジー群の合成積を用いた代数の構成という手法によ
る．これは幾何学的表現論で用いられてきた手法である．表現論と関連した
多様体の新しいクラスとして，クーロン枝は物理との関係を抜きにしても，純
粋に数学的な対象としてもおもしろいものである．

\cite{2016arXiv160403625B}において，箙ゲージ理論と呼ばれるゲージ
理論のクーロン枝の，さらにそのホモロジー(正確には交叉コホモロジー)を考
えることにより，Kac-Moodyリー環の表現が構成されると予想した．
これは，Kac-Moodyリー環が有限次元の複素単純リー環の場合には，クーロン枝
とアファイン・グラスマンとの関係\cite{2016arXiv160403625B}を使うと，
幾何学的佐武対応とよばれている構成に他ならない．幾何学的佐武対応が，
アファイン・リー環に拡張されることは Braverman-Finkelberg
\cite{braverman-2007}により予想されていた
が，\cite{2016arXiv160403625B}の予想は，これをさらにKac-Moodyリー環に拡
張するものである．この論説の執筆時点で，$A$型のアファイン・リー環の場合
に，予想の証明が与えられている \cite{2018arXiv181004293N}．
一方，この予想は，箙多様体のホモロジー群の上にKac-Moodyリー環の表現を構
成した筆者の研究~\cite{Na-quiver,Na-alg}と形式的に似ている．クーロン枝と箙多様
体(より一般的な状況では\subsecref{subsec:Higgs}で説明するヒッグス枝)と
の間には，シンプレクティック双対
性\cite{MR3594664,2014arXiv1407.0964B}と呼ばれる不思議な関係があること
が予想されており，Kac-Moodyリー環の表現が二つの多様体の上に与えられるこ
とは，その反映であると考えられている．

なお，物理的な背景の説明は省略するので，興味を持たれた方
は\cite{2016arXiv161209014N}や，その加筆英訳
版\cite{2017arXiv170605154N}をお読みいただきたい．

\section{クーロン枝とその量子化の定義}

この節では，\cite{2016arXiv160103586B}で与えられた，ゲージ理論に付随したクー
ロン枝の数学的に厳密な定義を紹介する．また，副産物として，その量子化も導
入する．

\subsection{幾何学的佐武対応の簡単な復習}\label{subsec:geom_satake}

定義のための準備に必要なことと，あとの\secref{sec:Coulomb_Kac_Moody}で
幾何学的佐武対応のKac-Moodyリー環への拡張を考えることから，
通常の複素簡約群 $\bG$ の場合について簡単に復習しておく．

幾何学的佐武対応に関する重要な結果とし
て \cite{Lus-ast,1995alg.geom.11007G,Beilinson-Drinfeld,MV2}があげられ
る．\secref{sec:Coulomb_Kac_Moody}では，Mirkovi\'c-Vilonen による定式
化\cite{MV2}を，一般化する．
幸いにして，\cite{MV2}についてたくさんの読みやすい解説があるので，ここ
では最小限の復習にとどめても，簡単に補うことができると思う．
なお，群ではなくリー環の幾何学的佐武対応と言っているのは，箙ゲージ理論
のクーロン枝が，ディンキン図形の組み合わせ論にのみ依存して決まっている
ためである．表現としては任意の可積分表現が現れるので，通常の場合でいう
と，$\bG$は随伴群，Langlands双対$\bG^\vee$は単連結なものを考えてい
る．

$\bT$を$\bG$の極大トーラスとする．以下の話では，$\bG = \GL_n(\CC)$,
$\bT$ は対角行列の全体であると思って読んでいただいて構わない．

$\cO = \CC[[z]]$,
$\cK = \CC((z))$を形式的級数のなす環，形式的ローラン級数のなす体とし，
対応する空間，形式的円盤 $D = \Spec \cO$と，形式的穴あき円盤 $D^\times = \Spec
\cK$をとる．
$\bG_\cK = \bG(\CC((z)))$, $\bG_\cO = \bG(\CC[[z]])$とする．

$\bG$のアファイン・グラスマン$\Gr_\bG$は，商空間$\Gr_\bG = \bG_\cK/\bG_\cO$
として定められる．簡約群 $\bG$ の旗多様体 $\bG/\mathbf P$ ($\mathbf
P$ は放物部分群) の無限次元リー群 $\bG_\cK$における類似，とみなすこと
ができるが，むしろ $\bG/\mathbf P$ の
$\mathbf P$軌道の閉包の，シューベルト多様体の類似と考えた方がより正確であ
る．
$\bG$の余ウェイト$\lambda$に対して，$\lambda(z)\in\bT_\cK$は，
$\Gr_\bG$の点を与える．これを$z^\lambda$と表す．
$\Gr_\bG$の$\bG_\cO$軌道は，支配的な余ウェイト$\lambda$によって
$\bG_\cO\cdot z^\lambda$と書くことができる．閉包の包含関係が誘導する順序は，
支配的順序になる．$\lambda$に対応する軌道 $\Gr_\bG^\lambda$ の閉
包$\overline{\Gr}_\bG^\lambda$は，一般には特異点を持つが有限次元の射影
代数多様体であり，以下の構成も，$\Gr_\bG$を直接扱うのではなく，厳密には$\overline{\Gr}_\bG^\lambda$で
議論しておいて，$\lambda$について極限をとるという形をとる．これ
が，$\Gr_\bG$はシューベルト多様体の類似と思ったほうがより正確といった理
由である．

$\bG$の極大コンパクト部分群 $\bG_c$ を取り，$S^1$から$\bG_c$への多項式
写像で$1\in S^1$を $\bG_c$の単位元にうつすもの全体を $\Omega\bG_c$ で表
わす．基点付きループ群などと呼ばれる．このと
き$\Omega\bG_c$と$\Gr_\bG$は同相であることが知られてい
る．(\cite[\S8.3]{MR900587}参照)
$1$を単位元に移すという条件を外したもの
を$\operatorname{Map}(S^1,\bG_c)$で表すと，$\Omega\bG_c =
\operatorname{Map}(S^1,\bG_c)/\bG_c$であり，これは，$\bG/\mathbf P
\cong \bG_c/(\bG_c\cap\mathbf P)$ の右辺の類似である．

幾何学的佐武対応は，$\Gr_{\bG}$の幾何を用いて，Langlands双対
群$\bG^\vee$の有限次元表現を構成するものである．ここでは，複素数体上の表
現を考えることにする．従って，有限次元表現は既約表現の直和に分解し，
既約表現の同型類は支配的なウェイトでパラメトライズされる．
Langlands双対群$\bG^\vee$の定義は省略するが，$\bG^\vee$の極大トーラ
ス$\bT^\vee$は$\bT$
の双対トーラスであり，$\bT$の余ウェイ
ト$\lambda\in\Hom(\CC^\times,\bT)$
は$\bT^\vee$のウェイト$\in\Hom(\bT^\vee,\CC^\times)$とも思えることは注
意しておこう．

余ウェイト$\lambda$に対応する$\bG^\vee$の既約表現を$V(\lambda)$と書く．
\cite{MV2}では，$V(\lambda)$のウェイト空間 $V_\mu(\lambda)$ を次のよう
に構成している．$\mu$を$\bG$の余ウェイトと思い，点$z^\mu\in\Gr_{\bG}$を
考える．単純余ルート$\alpha_i$を用いて$\lambda-\mu = \sum
v_i\alpha_i$と表したとき，$z^\mu\in\overline{\Gr}_{\bG}^\lambda$な
らば$v_i$は非負である．$\nu$を支配的で正則な余ウェイトとしたと
き，repelling set
\begin{equation}\label{eq:h:1}
  \left\{ x\in \overline{\Gr}_{\bG}^\lambda
    \,\middle|\,
    \lim_{t\to \infty} \nu(t)x = z^\mu
  \right\}
\end{equation}
は，Mirkovi\'c-Vilonenサイクルと呼ばれる既約成分を持ち，それらは$\sum
v_i$次元である．このときウェイト空間$V_\mu(\lambda)$は，repelling setの
最高次のホモロジー群として実現される．なお，本論説ではホモロジーといっ
たときは，Borel-Mooreホモロジー群のことと約束する．Mirkovi\'c-Vilonenサ
イクルはコンパクトではないが，その基本類はBorel-Mooreホモロジー群の元と
して定義される．

次にテンソル圏を用いた定式化を復習しよう．
$\Gr_\bG$の上の複素係数の$\bG_{\cO}$同変な偏屈層のなすアーベル
圏 $\Perv_{\bG_\cO}(\Gr_\bG)$ を考える．その対象とし
て，$\Gr_\bG^\lambda$の定数層から定まる交叉コホモロジー
層 $\IC(\overline{\Gr}_{\bG}^\lambda)$ がある．その台
は$\overline{\Gr}_{\bG}^\lambda$ である．$\Perv_{\bG_\cO}(\Gr_\bG)$ は，
支配的な余ウェイト$\lambda$ に対応す
る$\IC(\overline{\Gr}_{\bG}^\lambda)$を単純な対象として持つ，半単純なアー
ベル圏である．さらに，$\Perv_{\bG_\cO}(\Gr_\bG)$ は合成積によるテンソル
圏の構造が入ることが知られている．雑にいうと，$\Gr_\bG$ を基点付きルー
プ群$\Omega\bG_c$と思ったときの積写像$m\colon
\Gr_\bG\times\Gr_\bG\to\Gr_\bG$による押し出し準同型 $m_*$ により定義さ
れ
る．$\Perv_{\bG_\cO}(\Gr_\bG)\boxtimes\Perv_{\bG_\cO}(\Gr_\bG)$ が
$m_*$ で$\Perv_{\bG_\cO}(\Gr_\bG)$ に移されることは，非自明な主張であ
り，$m$ が semi-small という性質を持つことの帰結である．

このとき
$(\Perv_{\bG_\cO}(\Gr_\bG),m_*)$ と，$\bG$のLanglands双対
群 $\bG^\vee$ の有限次元複素表現の全体のな
す$(\operatorname{Rep}_{\bG^\vee},\otimes)$が，テンソル圏として同値であ
るというのが，幾何学的佐武対応の主たる主張である．ここで，$\otimes$は表
現のテンソル積である．さらに，圏同値のもと
で $\IC(\overline{\Gr}_{\bG}^\lambda)$ は $\lambda$ を最高ウェイトとす
る$\bG^\vee$の既約表現 $V(\lambda)$に対応する．


$\IC(\overline{\Gr}_{\bG}^\lambda)$と，Mirkovi\'c-Vilonenサイクルの基本
類との関係は，偏屈層の双曲制限関手というもので与えられ
る．\eqref{eq:h:1}の最高次のホモロジーの$\mu$に関する直和
に，$\bG^\vee$の表現の構造を入れるためには，偏屈層と双曲制限関手を用い
る必要がある．
この論説では，正確な定義を与えないので，以下では偏屈層は用いない．
%

\subsection{三つ組の多様体}\label{subsec:triples}

ここから，クーロン枝の定義に入る．

引き続き $\bG$ は複素簡約群とし，$\fg$ を $\bG$のリー環とする．
$\bG$の極大トーラス$\bT$のリー環を$\mathfrak t$と表す．また $\mathbf
W$をワイル群とする．
%



次に，$\bN$を$\bG$の有限次元表現とする．ただし $\bN$は既約でなくてもよ
く，$0$であってもよい．
$\bN$とその双対表現$\bN^*$の直和$\bM = \bN\oplus\bN^*$ は$\bG$のシンプ
レクティックな表現であり，最近の研究
\cite{bdfrt,2022arXiv220901088T}で組$(\bG,\bM)$に対してクーロン枝が
定義されることが分かっているが，以下を踏まえたものなので，
ここでは$(\bG,\bN)$に対して定義を与える．
\footnote{正確には
  障害類$\pi_4(\bG)\to\pi_1(\grpSp(\bM))$が消えるという条件が必要
  である．
}

$\bN$に値を持つ形式的ローラン級数，形式的級数$\bN((z))$, $\bN[[z]]$をそ
れぞれ $\bN_\cK$, $\bN_\cO$ で表す．
表現$\bN$に付随したアファイン・グラスマン上のベクトル
束$\cT$を，$\bG_\cK\times^{\bG_\cO}\bN_\cO$で定める．正確に
は，$s(z)\in\bN_\cO$の展開を途中で止め，上のよう
に$\overline{\Gr_\bG^\lambda}$に制限することによって，$\cT$は
射影多様体上のベクトル束の逆極限の直極限になる．
$\Pi\colon\cT\to\bN_\cK$ を
$[g(z),s(z)]\mapsto g(z)s(z)$で定義する．ここで，$g(z)\in\bG_\cK$,
$s(z)\in\bN_\cO$ であ
り，$[g(z),s(z)]$ は，$(g(z),s(z))$ の
$\bG_\cK\times^{\bG_\cO}\bN_\cO$における代表元であ
る．$s(z)$ は，$z=0$で極を持たないが，$g(z)$ は一般に持つので，掛けたも
の $g(z)s(z)$ は極を持つ可能性があり，そのために$\Pi$の像は $\bN_\cK$
になる．

$\cT$の閉部分多様体$\cR$として，$g(z)s(z)$が$z=0$で極を持たないという条件を課して，定められる空間と定義する．
\begin{equation*}
  \cR = \{ [g(z), s(z)]\in\cT \mid \Pi([g(z),s(z)]) = g(z)s(z) \in\bN_\cO \}.
\end{equation*}
この空間$\cR$は，ここでは説明しないモジュライ空間としての構成の理由によ
り，{\bf 三つ組の多様体(variety of triples)}と呼ばれる．

$\bG_\cO$が$\cR$に，
\begin{equation}\label{eq:1}
  \cR\ni [g(z), s(z)] \mapsto [h(z)g(z), s(z)]
  = [h(z) g(z) h(z)^{-1}, h(z) s(z)] \quad
  h(z)\in\bG_\cO
\end{equation}
により作用する．

この空間の定義の意味を少し説明するために，有限次元の状況で類似する空間を考えてみ
る．$\bG$を上と同様に複素簡約群，$\mathbf P$ をその放物部分
群，$\mathbf V$を$\bG$の表現，$\mathbf V'$ を$\mathbf P$で不変
な $\mathbf V$の部分空間とする．このとき
\begin{equation*}
  \cR = \{ [g,s]\in \bG\times^{\mathbf P}\mathbf V' \mid
  gs \in\mathbf V'\}
\end{equation*}
が，上の$\cR$の有限次元における類似である．(簡潔さのために同じ記号$\cR$で表
した．) $\mathbf V'$は，$\bG$で不変とは限らないから，$gs$ は一般に
は，$\mathbf V$ に入るだけで，$\mathbf V'$ には入らず，そのため
に$gs\in\mathbf V'$が非自明な条件になっている．
上では，$\bG$, $\mathbf P$, $\mathbf V$, $\mathbf V'$ をそれぞ
れ $\bG_\cK$, $\bG_\cO$, $\bN_\cK$, $\bN_\cO$ と取ったものになってい
る．
$\cT$の類似は，$\bG\times^{\mathbf P}\mathbf V'$である．
この状況のもとで，上の空間よりも大きな空間
\begin{equation*}
  \{ ([g_1,s_1], [g_2, s_2]\in
  \cT
  \times
  \cT
  \mid g_1s_1 = g_2 s_2 \}
\end{equation*}
を考えてみる．これは，$\cT
$ をそれ自身と$\mathbf V$ の上でファイバー積を取ったもの
\begin{equation*}
  \cT\times_{\mathbf V}\cT
\end{equation*}
に他ならない．ここで，$\cT
\to \mathbf V$は，上の
写像$\Pi$の類似で，$[g,s]\mapsto gs$で与えられるものに他ならない．

たとえば $\mathbf P$をBorel部分群$\mathbf B$ として，$\mathbf
V$ を随伴表現 $\mathfrak g$，$\mathbf V'$を上三角部分環$\mathfrak n$
と取ると，$\cT
$は，旗多様体 $\mathbf G /\mathbf B$の余接束であり，$\cT\times_{\mathbf V}\cT$
はワイル群のSpringer対応に現れるSteinberg多様体になる．

上の $\cR$は，$\cT\times_{\mathbf V}\cT$の部分多様体で，$g_2 =
\operatorname{id}$ となるものに他ならない．実際，$g_2 =
\operatorname{id}$ であれば $g_1 s_1 = s_2$で，$s_2$を消去して，その代
わりに $g_1s_1\in\mathbf V'$ を課せば，$\cR$の定義に他ならない．
$\cT\times_{\mathbf V}\cT
$には，$\bG$が作用しており，$\bG$同変な幾何学は部分空間 $\cR
$の$\mathbf P$同変な幾何学と同じである．象徴的に書けば，
\begin{equation}
  \left[
    \bG \backslash \left(
      \cT\times_{\mathbf V}\cT
    \right)\right] =
  \left[\mathbf P\backslash \cR
  \right]
\label{eq:2}\end{equation}
が成り立つ．

もともとの無限次元の状況では，$\cT\times_{\mathbf V}\cT$を直接に扱うこ
とは技術的に困難と思われるので，その代わりに$\cR$を考えている．しか
し，$\cT\times_{\mathbf V}\cT$が背景にあることをポイントとして押さえて
おくことは，以下の構成を理解する上では重要になる．

\begin{Remark}\label{rem:Hessenberg}
  \cite[\S5.4]{oy}では，有限次元の場合の$\cT\to \mathbf V$のファイバー
  をHessenberg多様体とよんでいる．(\cite{MR2209851}も参照．) これは，
  Springer fiberの一般化である．従って，三つ組の多様体
  からできる$\cT\to\bN_\cK$のファイバーは，無限次元リー群 $\bG_\cK$ に
  おけるHessenberg多様体の類似と考えることができる．
\end{Remark}

\subsection{合成積}

次に$\cR$の$\bG_\cO$同変ホモロジー群$H^{\bG_\cO}_*(\cR)$を考える．厳密
には，$\cT$の原点におけるファイバーの基本類が次数$0$になるように，次数
をうまく定義する必要があるが，詳細は略す．また，奇数次のホモロジーが消
えていること，$H^*_\bG(\mathrm{pt})$上自由な加群になっていることなどは，
アファイン・グラスマン多様体のシューベルト胞体分割を考えると，ただちに
従う．

有限次元の類似においては，同変ホモロジー群の誘導により，\eqref{eq:2}から
\begin{equation}\label{eq:3}
  H^{\mathbf P}_*(\cR) \cong H^{\bG}_*(\cT\times_{\mathbf V}\cT)
\end{equation}
となることに注意しよう．このとき，$(i,j)$成分への射影
\begin{equation*}
  \cT\times_{\mathbf V}\cT\times_{\mathbf V}\cT \xrightarrow{p_{ij}}
  \cT\times_{\mathbf V}\cT \qquad (i,j) = (1,2), (2,3), (1,3)
\end{equation*}
を用いて
\begin{equation*}
  c\ast c' = p_{13*}(p_{12}^*c \cap p_{23}^*c')
\end{equation*}
と $H^{\mathbf G}_*(\cT\times_{\mathbf V}\cT)$上に合成積が定義される．
ここで，$\cap$ は $\cT\times\cT$における(台を考慮した)キャップ積であ
り，$\cT$が滑らかな多様体であることを用いて定められている．

我々の無限次元の状況では，$H^{\bG_\cK}_*(\cT\times_{\bN_\cK}\cT)$ が，
定義されていないこと，$\cT$が滑らかではないことにより，上の定義は適用できないが，
これを回避して $H^{\bG_\cO}_*(\cR)$ に直接，合成積
\begin{equation*}
  \ast\colon H^{\bG_\cO}_*(\cR)\otimes H^{\bG_\cO}_*(\cR)
  \to H^{\bG_\cO}_*(\cR)
\end{equation*}
が定義される．厳密な定義は，技術的なのでここでは略す．
$\bN = 0$ のときは，\cite{MV2}が幾何学的佐武対応の構成に用いた図式を，
同変ホモロジー群に適用したものに他ならない．この特別な場合の構成
は，\cite{MR2135527} で行われており，また，得られる環の具体的な記述も与
えられている．
なお，\cite{MR2135527}では随伴表現$\bN=\fg$の場合も，$K$理論版のときに，
取り扱われている．

合成積代数
$H^{\bG}_*(\cT\times_{\mathbf V}\cT)$の単位元は，対角線集
合$\Delta_\cT$の基本類であるが，\eqref{eq:3}の同型のも
と，$\cR$の$\bG/\mathbf P$の原点のファイバー
$\{ [g,s]\in \bG\times^{\mathbf P}\mathbf V'\mid g \in\mathbf P\}$ の基
本類に対応することに注意しよう．$H^{\bG_\cO}_*(\cR)$においても，同様
に $\Gr_{\bG} = \bG_\cK/\bG_\cO$ の $\bG_\cO$のファイバーの基本類が単位
元になっている．

このとき次が成立する．
\begin{Theorem}
  $(H^{\bG_\cO}_*(\cR),\ast)$ は可換環である．
\end{Theorem}

合成積で環を構成する手法は，幾何学的表現論で広く使われており，たとえば
途中で紹介したSteinberg多様体の場合には，(退化した)アファイン・ヘッ
ケ環が実現される．
%
この例では，得られる環は，非可換環であり，合成積の一般論からは$\ast$が
可換になる理由はなく，上の定理は今の状況の特殊性を表している．

幾何学的佐武対応を思い起こせば，可換性は不思議ではない．幾何学的佐武対
応で
は，$(\operatorname{Perv}_{\bG_\cO}(\Gr_{\bG}),\allowbreak m_*)$ と
$(\operatorname{Rep}_{\bG^\vee},\otimes)$ がテンソル圏同値となったが，
後者のテンソル圏は可換，すなわち $V\otimes W\cong W\otimes V$ であるの
で，前者もそうである．この同型を幾何学的に説明するの
がBeilinson-Drinfeldによるアファイン・グラスマンの1パラメータ変形であ
り，同じアイデアを使って上の定理が証明され
る．\cite{2017arXiv170602112B}を参照．(\cite{2016arXiv160103586B}では，
計算による直接証明も与えている．)

さて，$(H^{\bG_\cO}_*(\cR),\ast)$は可換環になったので，そのスペクトラム
としてアファイン多様体を導入することができる．これが，クーロン枝の数学
的な定義である．
\begin{equation*}
   \cM_C = \Spec (H^{\bG_\cO}_*(\cR),\ast)
\end{equation*}
さらに，$(H^{\bG_\cO}_*(\cR),\ast)$が有限生成であることや，整であること
を証明できるので，$\cM_C$は既約なアファイン多様体である．また正規である
ことも示されている．

\subsection{ヒッグス枝}\label{subsec:Higgs}

  シンプレクティック・ベクトル空間 $\bM=\bN\oplus \bN^*$ への $\bG$作用
  の運動量写像$\mu\colon \bM\to\fg^*$ は
  \begin{equation*}
    \langle \mu(x,y),\xi\rangle = \langle y, \xi x\rangle
  \end{equation*}
  で与えられる．このとき，{\bf ハミルトニアン簡約}
  $\bM\tslash \bG$
  は，$\mu^{-1}(0) \dslash \bG$で定義される．ここで，$\dslash$はカテゴ
  リカル商 $\Spec \CC[\mu^{-1}(0)]^\bG$ である．$\fg^*$ の余随伴作用で
  固定される元$\zeta$におけるレベルに取り替え
  た $\mu^{-1}(\zeta)\dslash \bG$や，$\bG$の乗法的指
  標 $\chi\colon\bG\to\BG_m$を取り，幾何学的不変式論による
  商
  $\mu^{-1}(0)\dslash_\chi \bG = \Proj\left( \bigoplus_{n\ge 0}
    \CC[\mu^{-1}(0)]^{\bG,\chi^n}\right)$ を考えることもよくある．ここ
  で$[\ ]^{\bG,\chi^n}$ は相対不変式を表す．
  物理のゲージ理論の文脈では，これらは{\bf ヒッグス枝}と呼ばれる．なお，
  ハミルトニアン簡約は$\bM$に対して定義され，必ずしも$\bN\oplus\bN^*$と
  分解している必要はないことを注意しておく．

  共通の$\bG$, $\bN$ に対して定義されるクーロン枝とヒッグス枝の間には一
  見すると何も関係がないように思われる．実は両者の間には不思議な関係が
  あることが見つけられており，{\bf シンプレクティック双対性}と名づけられている．
  \cite{MR3594664,2014arXiv1407.0964B}参照．

\subsection{量子化されたクーロン枝}\label{subsec:quantum}

おもしろいことに，$H^{\bG_\cO}_*(\cR)$は，その非可換変形を同時に作るこ
とができる．実際，形式的円盤 $D$に，$\CC^\times$が
$z\mapsto tz$により作用する．この作用は，ループ・ローテーションと呼ばれ
るが，今まで使ってきた様々な空間への作用を引き起こす．特に，$\bG_\cO$に
作用して，半直積$\bG_\cO\rtimes\CC^\times$を考えることがで
き，$\cR$に$\bG_\cO\rtimes\CC^\times$が作用する．そこで，同変ホモロジー
群$H^{\bG_\cO\rtimes\CC^\times}_*(\cR)$を考え，合成積を同じように導入する．
こうして{\bf 量子化されたクーロン枝}を
\begin{equation*}
   \cAh = (H^{\bG_\cO\rtimes\CC^\times}_*(\cR),\ast)
\end{equation*}
と定義する．これは，非可環にな
り，$H^{\bG_\cO}_*(\cR)$ の $H_{\CC^\times}^*(\mathrm{pt}) =
\CC[\hbar]$ でパラメトライズされた非可環変形になる．

非可環変形が与えられると
\begin{equation*}
  \{ f, g\} =
  \left.\frac{\tilde f \ast \tilde g - \tilde g\ast \tilde f}\hbar
    \right|_{\hbar = 0}
\end{equation*}
により，$H^{\bG_\cO}_*(\cR)$はポアソン代数になる．ここで，$\tilde f$,
$\tilde g$ は$H^{\bG_\cO\rtimes\CC^\times}_*(\cR)$への持ち上げである．さ
らに，これは $\cM_C$の非特異部分上にシンプレクティック形式を定める．この意味で，
$\cAh$はクーロン枝$\cM_C$の{\bf 非可換変形}，すなわち量子化である．

先にあげた\cite{MR2135527}の続編の\cite{MR2422266}では，$\bN=0$の場合に
量子化を調べ，$\bG$のLanglands双対の戸田格子と同定しているが，詳細は
略す．
また，\cite{MR3013034}では，アファイン・グラスマンの代わりにアファイン
旗多様体，同変ホモロジー群の代わりに同変$K$群が用いられているが，$\bN =
\mathfrak g$の場合を扱っている．出てくる代数は，Cherednikの二重アファイ
ン・ヘッケ代数(DAHA)である．クーロン枝のようにアファイン・グラスマン多
様体にすれば，そのspherical partになり，同変ホモロジー群になれば楕円版
の代わりに三角関数版のDAHAになる．対応するクーロン枝
は$\ft\times\bT^\vee/\mathbf W$であり，あとで定理~\ref{thm:classical}で
述べるように，これは量子補正を受ける前の古典的なクーロン枝と同じである．
これらの結果は，クーロン枝の定義の先駆である．

\subsection{可積分系}

$H^{\bG_\cO}_*(\cR)$ は，同変ホモロジー群であり，点の同変コホモロジー
環 $H^*_{\bG}(\mathrm{pt})$からの環準同型を持つ．よって，
\begin{equation*}
  \varpi\colon \cM_C \to \Spec H^*_{\bG}(\mathrm{pt})
\end{equation*}
という射が誘導される．よく知られている事実により，右辺は随伴
商 $\mathfrak g\dslash \bG$ である．また，随伴商は
$\fg$のカルタン部分環をワイル群 $\mathbf W$で割った
$\mathfrak t/\mathbf W$ と同型であることもよく知られている．

環準同型 $H^*_{\bG}(\mathrm{pt}) \to H^{\bG_\cO}_*(\cR)$ は，量子化
に$H^*_{\bG\times\CC^\times}(\mathrm{pt}) \to
H^{\bG_\cO\rtimes\CC^\times}_*(\cR)$と拡張される．これは単射であること
が示され，したがっ
て$H^*_{\bG\times\CC^\times}(\mathrm{pt})$は，$\cAh$の可換部分環になる．
よって，$H^*_{\bG}(\mathrm{pt})$はポアソン可換な$H^{\bG_\cO}_*(\cR)$の部
分環である．すなわち，$H^*_{\bG}(\mathrm{pt})$の二つの元のポアソン括弧
は$0$になる．

なお，$H^*_{\bG\times\CC^\times}(\mathrm{pt})$は可換部分環にはなるが，
中心ではない．有限次元の類似の同型 \eqref{eq:3}におい
て $H^*_{\bG}(\mathrm{pt})$が中心に含まれることは，合成積の定義に現れる
射影 $p_{ij}$ が $\bG$同変であることから従うが，
有限次元の類似の$H^*_{\bG}(\mathrm{pt})$に，今の状況で対応させるべきも
のは$H^*_{\bG_\cK\rtimes\CC^\times}(\mathrm{pt})$ であり，これ
は$H^*_{\bG_\cO\rtimes\CC^\times}(\mathrm{pt})$とは異なるし，そもそもど
のように定義したらいいのかさえも明らかでない．これは，形式的な説明で
あるが，\subsecref{subsec:torus}で見るように，簡単な具体例でも中心には入らな
いことが計算でチェックできる．

\subsection{トーラスの作用}\label{subsec:hamtori}

合成積は，ホモロジーの次数付けと整合的であ
り，$H^{\bG_\cO\rtimes\CC^\times}_*(\cR)$は$\ZZ$で次数付けられた環にな
る．また，この次数は，$\cM_C = \Spec
H^{\bG_\cO}_*(\cR)$への$\CC^\times$作用を定める．実際，次数$d$の部分空間
が，$\CC^\times$作用に関してウェイト分解したときの次数$d$の部分空間にな
るように，$\CC^\times$作用が定まる．

さらに，$\cR$の連結成分の集合$\pi_0(\cR)$ は
\begin{equation*}
  \pi_0(\cR) = \pi_0(\Gr_{\bG}) = \pi_0(\Omega\bG_c) \cong \pi_1(\bG)
\end{equation*}
と同一視され，ホモロジー群は
\(
   H^{\bG_\cO\rtimes\CC^\times}_*(\cR)
   = \bigoplus_{\gamma\in\pi_1(\bG)} H^{\bG_\cO\rtimes\CC^\times}_*(\cR_\gamma)
\)
と分解し，$\pi_1(\bG)$で次数付けられた環になる．$\cM_C$に作用するの
は，$\pi_1(\bG)$のポントリャーギン双対
$\pi_1(\bG)^\wedge$となる．$\bG$が半単純のときには，$\pi_1(\bG)$は有限
群になることがよく知られているので，$\pi_1(\bG)^\wedge$も有限群となって
しまい，あまりおもしろくないが，$\bG$が一般線形群 $\GL(n)$のときに
は$\pi_1(\bG)^\wedge = \ZZ^\wedge = \CC^\times$ となり，い
ろいろな役割を果たすことになっていく．さらに，\secref{sec:quiver}で取り扱う箙ゲージ理
論のクーロン枝の場合には，$\bG$は一般線形群の直積であ
り，$\pi_1(\bG)^\wedge$はトーラスになる．これは，箙が定める群の極大トー
ラスと同一視することができる．

\subsecref{subsec:Higgs}で言及したシンプレクティック双対性の枠組みにお
いての$\pi_1(\bG)^\wedge$の役割を紹介しよう．$\bG$が一般線形群の直積の
場合，$\Hom_{\mathrm{grp}}(\bG,\BG_m) \cong
\Hom_{\mathrm{grp}}(\BG_m,\pi_1(\bG)^\wedge)$
が$\pi_1(\BG_m)^\wedge = \ZZ^\wedge = \BG_m$により誘導される．したがっ
て，$\pi_1(\bG)^\wedge$の余ウェイトは，$\bG$の乗法的指標と思うことができる．
前者は，クーロン枝に作用する1パラメータ変換群と思うことができる．一方，
\subsecref{subsec:Higgs}で述べたように，後者はヒッグス枝の変形や幾何学的不変
式論の商を与える．これらは，幾何学的にはまったく性質の異なるものである
が，シンプレクティック双対性を定式化する際には，両者を関係させることに
なる．

\subsection{フレーバー対称性}\label{subsec:flavor}

Steinberg多様体の合成積代数の場合には，アファイン・ヘッケ環が実現
されると述べた．正確には，上で考えていた$H^{\bG}_*(\cT\times_{\mathbf
  V}\cT)$では，実現されるのはアファイン・ワイル群の群環であって，アファイ
ン・ヘッケ環を実現するためには，$\cT = \bG\times^{\mathbf P}\mathbf V'$のファイバー
に作用する$\CC^\times$作用を付け加えて，$\bG\times\CC^\times$同変ホモ
ロジー群を考える必要がある．この作用は，ループ・ローテーションからく
る $\CC^\times$作用とは異なることは注意しておく．

同様の構成は，より一般の状況で考えることができ，物理の文脈でフレーバー
対称性といわれている．はじめに$\bN$は$\bG$の表現であるとしたが，$\bG$を
正規部分群として含むような別の複素簡約群$\tilde\bG$の表現の制限になって
いると仮定しよう．このとき，\eqref{eq:1}の右辺の第二項
で $h(z)\in\tilde\bG_\cO$としても作用が well-defined であることに注意し
て，大きな群に関する同変ホモロジー
群$H^{\tilde\bG_\cO}_*(\cR)$を考えることができる．これ
は，$H^{\bG_\cO}_*(\cR)$の
$H_{\tilde\bG/\bG}^*(\mathrm{pt})$でパラメトライズされた変形を与える．

また，群とその表現をはじめから，$\tilde\bG$と$\bN$としたものを考え
る．三つ組の多様体も代わるので，$\cR_{\tilde\bG,\bN}$ で表す．
さらに，$\tilde\bG/\bG$はトーラスであると仮定し，$\bT_F$ で表
す．($F$ はフレーバーを意味する．) すると，$\Spec
H^{\tilde\bG}_*(\cR_{\tilde\bG,\bN})$に
\begin{equation*}
  \bT_F^\vee = \pi_1(\bT_F)^\wedge
  \to \pi_1(\tilde\bG)^\wedge
\end{equation*}
を通じて，双対トーラス $\bT_F^\vee$ が作用する．このとき次が成り立つ．

\begin{Proposition}\label{prop:reduction}
  $\cM_C = \Spec H^{\bG}_*(\cR)$
  は，$\Spec H^{\tilde\bG}_*(\cR_{\tilde\bG,\bN})$ の $\mathbf
  T_F^\vee$ によるハミルトニアン簡約である．
\end{Proposition}

$\Spec H^{\tilde\bG}_*(\cR_{\tilde\bG,\bN})$ から $\Spec
H^{\bG}_*(\cR)$ を得るためには，1) 同変ホモロジーをとる群
を $\tilde\bG$ から $\bG$ に取り替える，2) $\cR_{\tilde\bG,\bN}$ の連結
成分のうちで $\pi_1(\tilde\bG)\to\pi_1(\bT_F)$ の核に入っているも
のに対応する部分だけをとる，という二段階を経る必要がある．前者は，運動
量写像の値が$0$の部分多様体を取ることに対応し，後者は $\mathbf
T_F^\vee$ に関して，カテゴリカル商，すなわち，$\bT_F^\vee$不変部
分関数環の$\Spec$ をとることに対応する．

上で述べた変形$H^{\tilde\bG_\cO}_*(\cR)$は，運動量写像のレベルの値
を$0$から動かすことに対応する．また，$\bT_F^\vee$によるカテゴリカル商の
代わりに幾何学的不変式論における商，すなわち$\bT_F^\vee$の指標 $\nu$ を
とり，
\begin{equation*}
  \cM_C^\nu :=
  \Proj \bigoplus_{n\ge 0} H^{\bG}_*(\cR_{\tilde\bG,\bN})^{\bT_F^\vee,\nu^n}
\end{equation*}
と置き換えることができる．ここで添字の$\mathbf
T_F^\vee,\nu^n$は，$\bT_F^\vee$の作用で$\nu^n$倍される相対不変部分空間を意
味する．幾何学的不変式論の一般論から，射影的射 $\cM_C^\nu\to\cM_C$ が誘
導される．これは，$\cM_C$ の部分的特異点解消を与えることがチェックされ
ている．

\subsecref{subsec:Higgs}において，$\Hom_{\mathrm{grp}}(\BG_m,\bT_F)$ の元は，ヒッ
グス枝に働く1パラメータ変換群と思えることに注意しよう．実際，ヒッグス枝
は$\bG$による商空間なので，$\bT_F=\tilde\bG/\bG$が作用する．一方，上で
説明したのは，$\Hom_{\mathrm{grp}}(\BG_m,\bT_F)\cong
\Hom_{\mathrm{grp}}(\bT_F^\vee,\BG_m)$ がクーロン枝の変形や幾何学的不変
式論の商を与えることであった．これは，\subsecref{subsec:hamtori}の最後に述べ
た $\pi_1(\bG)^\wedge$の果たす役割の，ヒッグスとクーロン枝を入れ替えたも
のに他ならない．
物理の文脈では，FIパラメータ($\in\Hom_{\mathrm{grp}}(\bG,\BG_m)$)と質量
パラメータ($\in\Hom_{\mathrm{grp}}(\BG_m,\bT_F)$)の役割が，クーロン枝と
ヒッグス枝で入れ替わると説明されていたが，クーロン枝の数学的な定義では
この性質が数学的に厳密に確立されたわけである．

\section{群がトーラスの場合の例}

前節の構成は，無限次元空間のホモロジーを使うもので，抽象的な構成に見え
るかもしれないので，$\bG$がトーラスの場合に具体的に計算してみよう．

\subsection{トーラスの$0$表現の場合}\label{subsec:torus}

$\bG = \CC^\times$とし，$\bN = 0$ とする．これは，一番自明な例であ
る．$\bN=0$なので，$\cR$はアファイン・グラスマン$\Gr_\bG$に他ならない．
また$\Gr_{\CC^\times}$は，整数$\ZZ$でパラメトライズされた離散的な空間にな
る．実際，$g(z) = z^n$ ($n\in\ZZ$)が対応する点を表す．よって
\begin{equation*}
   H^{\bG_\cO}_*(\cR) = \bigoplus_n H^{\CC^\times}_*(\mathrm{pt})
\end{equation*}
となる．$H^{\CC^\times}_*(\mathrm{pt})\cong H^*_{\CC^\times}(\mathrm{pt})$は，
一変数の多項式環$\CC[w]$であ
る．これが各整数$n$の上に乗っているので，$m$の上の多項式と$n$の上の多項
式を掛けるとどうなるかを，合成積の定義に戻って計算する．
表現 $\bN=0$ のときは，合成積は幾何学的佐武対応のときと同様に群の掛け算
作用から決まるものを使って良い．従って
\begin{equation*}
  \Gr_{\CC^\times}\times\Gr_{\CC^\times}\to \Gr_{\CC^\times};
  z^m \times z^n \mapsto z^{m+n}
\end{equation*}
が，ホモロジー群に引き起こす押し出し準同型が$\ast$に他なら
ない．すると$m$の上の$f(w)$と$n$の上の$g(w)$を掛けたものは，$m+n$の上
の$f(w)g(w)$になる．すなわち，$n=1$ の上の$1$(基本類に対応する)を$x$で
表すと，
\begin{equation*}
   H^{\bG_\cO}_*(\cR) \cong \CC[w,x^{\pm 1}] = \CC[\BA\times\BA^\times]
\end{equation*}
となる．従って，今の場合のクーロン枝は$\BA\times\BA^\times$である．

もう一歩，精密に見るために$\bG$はトーラス$\bT$で，表現はやはり$0$であると
する．$\Gr_\bT$は離散的な空間で，$\Hom(\CC^\times,\bT)$でパラメトライズされ
ている．従って，$H^{\bT_\cO}_*(\cR) =
\bigoplus_{\lambda\in\Hom(\CC^\times,\bT)} H^*_\bT(\mathrm{pt})$であ
る．$H^*_\bT(\mathrm{pt})$は，$\bT$のLie環$\ft$上の多項式環$\CC[\ft]$である．
一方，$\lambda$に対応する元を$e^\lambda$と書くと，上と同様
に$e^\lambda\ast e^\mu = e^{\lambda+\mu}$ となる．これは，$\bT$の双
対 $\bT^\vee$ の指標 ($\Hom(\bT^\vee,\CC^\times) = \Hom(\CC^\times, \bT)$) と
見なすことができるから，クーロン枝は$\ft\times \bT^\vee = T^* \bT^\vee$であ
る．

量子化は，$w\in\ft^*$ と $e^\lambda\in \Hom(\CC^\times,\bT)$ で生成さ
れ，関係式は $w w' = w' w$, $e^\lambda e^\mu = e^{\lambda+\mu}$ と
\begin{equation*}
  [e^\lambda, w ] = e^\lambda w - w e^\lambda = \hbar \langle \lambda,w\rangle
  e^\lambda
\end{equation*}
である．$e^\lambda$ を$\ft$上の$\hbar$差分作用素
$f\mapsto f(\bullet+\lambda\hbar)$と見ることができるの
で，$\cAh$は，$\ft$上の多項式係数$\hbar$差分作用素のなす環に他ならな
い．
上の関係式は，合成積の定義に従って容易にチェックできるが，上の説明とは
異なり，単純
に$\Gr_{\CC^\times}\times\Gr_{\CC^\times}\to\Gr_{\CC^\times}$ から導か
れる写像ではなく，説明をサボっている定義を原論文に従って使う必要がある．

\subsection{$A$型単純特異点}\label{subsec:C2}

次に$\bG$は$\CC^\times$のままで, 表現を$\bN = \CC$ と標準表現に取ろう．
$\Gr_{\CC^\times}$は上で説明したように$\ZZ$でパラメトライズされる離散的な空間であり，$\cR$は各整数$n$の上にベクトル空間が乗っているものである．条件は$g(z) = z^n$によって原点に特異点が生じないというものであるから，
\begin{equation*}
  \cR = \bigsqcup_{n\in\ZZ} z^{n}\CC[z]\cap \CC[z] =
  \bigsqcup_{n\in\ZZ} z^{\max(0,n)}\CC[z]
\end{equation*}
である．各整数の上に乗っているものはベクトル空間であり，Thom同型により$H^{\bG_\cO}_*(\cR) \cong\allowbreak \bigoplus_n \allowbreak H_{\CC^\times}^*(\mathrm{pt})$となる．すなわちベクトル空間としては，\subsecref{subsec:torus}の例と同じである．しかし合成積は，\subsecref{subsec:torus}の例とは$n>0$の上のホモロジー類と$n<0$の上のホモロジー類の積が変わってくる．定義を省略したので，最後の結果だけいうと，$n=1$の基本類と$n=-1$の基本類を掛けたものが
\begin{equation*}
  z\CC[z] \to \CC[z]
\end{equation*}
の押し出し写像による，基本類の像になる．これは余次元$1$の部分空間である
から，同変ホモロジー群の元としては，$w$を基本類に掛けたものになる．従っ
て, $n=1$の基本類を$x$, $n=-1$の基本類を$y$とすると，$xy = w$が成り立つ．
この計算から
\begin{equation*}
  H^{\bG_\cO}_*(\cR) \cong \CC[w,x,y]/(w=xy) \cong \CC[x,y] = \CC[\BA^2]
\end{equation*}
が従う．よって今の場合のクーロン枝は$\BA^2$である．

表現をウェイト$1$の一次元表現の$\ell$個の直和に取り替えると，最後の部分の
計算が$z^{\ell}\CC[z]\to \CC[z]$のpushforwardに置き換わり，座標環
は$\CC[w,x,y]/(w^{\ell}=xy)$となる．これは，$A_{\ell-1}$型の単純特異点に他な
らない．

量子化については，\subsecref{subsec:torus}に帰着させるのが分かりやす
い．\subsecref{subsec:torus}の $n=1$の基本類を $r^1$, $n=-1$の基本類
を$r^{-1}$ とする．$r^1 r^{-1} = r^{-1} r^1 = 1$であった．このとき $x
= w^\ell r^1$, $y = r^{-1}$ によ
り，$\cAh$ は\subsecref{subsec:torus}の $\cAh$ の部分環として実現される．
従って，$xy = w^\ell$, $y x =
\begin{NB}
  r^{-1} w^\ell r^1 = r^{-1}(f \mapsto w^\ell f(w+\hbar))
  = f\mapsto (w-\hbar)^\ell f(w) =
\end{NB}%
(w-\hbar)^\ell$となる．ポアソン括弧は，$\{ x, y\} = \ell w^{\ell-1}$である．
\begin{NB}
  Since $x = z_1^\ell$, $y = z_2^\ell$, we have
  $\{ x, y\} = \{ z_1^{\ell}, z_2^\ell\} = \ell \{ z_1^\ell, z_2\} z_2^{\ell-1}
  = \ell^2 \{z_1, z_2\} (z_1 z_2)^{\ell-1}$
\end{NB}%

\subsection{トーリック超ケーラー多様体}

$\bG$をトーラス$\bT$とし，次の短完全列があるとする．
\begin{equation*}
  1 \to \bT \to (\CC^\times)^n \to \bT_F \to 1
\end{equation*}
ただし，$\bT_F$はやはりトーラスであると仮定する．このと
き，$(\CC^\times)^n$の標準表現$\CC^n$ を取り，$\bT$に制限して$\bN$とお
く．双対トーラスについての短完全列
\begin{equation*}
  1 \to \bT_F^\vee \to (\CC^\times)^n \to \bT^\vee \to 1
\end{equation*}
により，$\CC^n$は$\bT_F^\vee$の表現となることに注意しよう．

このとき
\begin{Theorem}[\protect{\cite[\S4(vii)]{2016arXiv160103586B}}]
  $(\bT, \bN)$ の定めるクーロン枝は，$\CC^n\oplus
  (\CC^n)^*$ の$\bT_F^\vee$によるハミルトニアン簡約である．
\end{Theorem}

\begin{proof}
  命題~\ref{prop:reduction}により，$(\bT,\bN)$の定めるクーロン枝は，
  $((\CC^\times)^n, \CC^n)$ の定めるクーロ
  ン枝の$\bT_F^\vee$によるハミルトニアン簡約である．\subsecref{subsec:C2}によ
  り，後者のクーロン枝は$\CC^n\oplus(\CC^n)^*$である．さら
  に，$\bT_F^\vee$の作用は，$((\CC^\times)^n)^\vee \cong
  (\CC^\times)^n$ の作用の制限であるが，\subsecref{subsec:C2}の計算か
  ら，$(\CC^\times)^n$の作用は標準的なものであることも従う．
\end{proof}

また，\subsecref{subsec:Higgs}に注意すると，この主張は$(\bT_F^\vee, \CC^n)$に付
随するヒッグス枝が，$(\bT,\bN)$のクーロン枝に等しい，と言っているとも理
解できる．

このような $\CC^n\oplus (\CC^n)^*$ をトーラスでハミルトニアン簡約して得
られる多様体は，{\bf トーリック超ケーラー多様体}と言われている．箙多様
体にもなっている特別な場合に \cite{MR1187554}で調べられ，一般の場合
に \cite{MR1792372}により導入された．
なお，先行研究を無視し，ハイパートーリック多様体と，\emph{間違って}呼ばれる
ことが多いので注意が必要である．

\section{局所化定理と古典的なクーロン枝}\label{sec:local}

\subsection{トーラス固定点集合}

次に同変ホモロジー群の局所化定理を用いて$H^{\bG_\cO}_*(\cR)$の解析を行
なおう．まず，同変ホモロジー群の一般的な性質から$H^{\bG_\cO}_*(\cR)
=H^{\bT_\cO}_*(\cR)^{\mathbf W}$に注意しよう．ただし，$\mathbf T$の正規
化群を $N(\bT)$ として，$H^{\mathbf
  T_\cO}_*(\cR)$にワイル群$\mathbf W = N(\bT)/\bT$が作用する．
さらに，$H^*_{\bT}(\mathrm{pt})$ の商体を
$\operatorname{Frac}H^*_{\bT}(\mathrm{pt})$とおくと，
同変ホモロジー群の局所化定理により
\begin{equation}\label{eq:4}
    H^{\bT_\cO}_*(\cR)\otimes_{H^*_{\bT}(\mathrm{pt})}
    \operatorname{Frac}H^*_{\bT}(\mathrm{pt}) \cong
  H^{\bT_\cO}_*(\cR^\bT)\otimes_{H^*_\bT(\mathrm{pt})}
  \operatorname{Frac}H^*_{\bT}(\mathrm{pt})
\end{equation}
となる．ただし，$\cR^\bT$は固定点集合で，
右辺から左辺への写像は，埋め込み写像
$\cR^\bT\hookrightarrow \cR$が誘導するものである．

従って，固定点集合を決定することが重要である．

\begin{Lemma}
    \textup{(1)} $\Gr_\bG$ の $\bT$作用に関する固定点集合は，$\Gr_\bT$ である．

    \textup{(2)} $\cR$の$\bT$作用に関する固定点集合は，群 $\bT$と表
    現 $\bN^\bT$に対する三つ組の多様体であり，さらに$\Gr_\bT\times
    (\bN^\bT)_\cO$ である．
\end{Lemma}

\begin{proof}
    (1)
    $\bT$の極大コンパクト部分群を$\bT_c$とすると，$\bT$固定点集合
    と，$\bT_c$固定点集合は同じである．さらに，同相 $\Gr_\bG \cong
    \Omega\bG_c$
    において，$\bT_c$ の作用は$\Omega\bG_c\ni c(z) \mapsto t c(z)
    t^{-1}$
    で与えられる．($z\in S^1$, $t\in\bT$) 従って，固定点集合
    は，$S^1$から$\bT$への多項式写像で，$1$を単位元に移すものの全
    体，$\Omega\bT_c$に他ならない．同
    相
    $\Gr_\bT\cong\Omega\bT_c$は，$\Gr_\bG\cong\Omega\bG_c$とcompatible
    で，固定点集合は$\Gr_\bT$である．

    (2) $\bN_\cO$, $\bN_\cK$ の$\bT$固定点集合は，それぞ
    れ$(\bN^\bT)_\cO$, $(\bN^\bT)_\cK$ である．よって(1)と合わせて結論
    の前半を得る．また，$\bT$の$\bN^\bT$への作用は自明であるから，後半
    を得る．
\end{proof}

$\bN^\bT$は$\bT$の自明表現で，クーロン枝には寄与せず，$\bN=0$の場合
の\subsecref{subsec:torus}と変わらない．これを，\eqref{eq:4}と合わせると

\begin{Theorem}[\protect{\cite[5.21]{2016arXiv160103586B}}]\label{thm:classical}
  埋め込み写像 $\cR^\bT\hookrightarrow\cR$は，双有理写像
  \begin{equation}\label{eq:5}
    \cM_C \approx (\ft\times \bT^\vee)/\mathbf W
  \end{equation}
  を誘導する．
\end{Theorem}

興味深いことに，物理では\eqref{eq:5}の右辺は古典的なクーロン枝として現
れ，$\cM_C$はその量子補正で得られる，と説明される．数学的な定義では，量
子補正されたものが始めから与えられる．

また，右辺は表現$\bN$にはよらず，$\bG$だけで決まっていることに注意しよう．

量子化されたクーロン枝$\cAh$の場合にも，同様に同変ホモロジーの局所化定
理を適用することができる．\subsecref{subsec:torus}において，$\ft\times
\bT^\vee$の量子化として$\ft$上の差分作用素の環が現われたことを思い出そ
う．これを局所化した有理関数係数の差分作用素の環に，$\cAh$が部分代数と
して実現されることが従う．
(\cite[Remark~5.23]{2016arXiv160103586B}参照．)

\section{箙ゲージ理論とそのクーロン枝}\label{sec:quiver}

\subsection{一般化された横断切片}\label{subsec:gen_slice}

現在のところ，クーロン枝の解析が一番進んでいるのは，$\bG$がトーラスの場
合を除くと，箙ゲージ理論の場合である．これは，箙 $Q =
(Q_0,Q_1)$ と，$Q_0$-graded な複素ベクトル空間の組$V = \bigoplus_{i\in
  Q_0} V_i$, $W = \bigoplus_{i\in Q_0} W_i$ に対して，
\begin{equation}\label{eq:7}
  \bG = \prod_{i\in Q_0} \GL(V_i), \quad
  \bN = \bigoplus_{h\in Q_1} \Hom(V_{\vout{h}}, V_{\vin{h}})
  \oplus\bigoplus_{i\in Q_0} \Hom(W_i, V_i)
\end{equation}
と取るものである．ここで，$h\in Q_1$は箙の辺で，$\vout{h}$ はその出発
点，$\vin{h}$ は到着点を表す．$\bG$は$\bN$に共役で作用
し，$\bN$は$\bG$の表現と見ることができる．

このとき，\subsecref{subsec:Higgs}により$(\bG,\bN)$に付随するヒッグス枝
は，\cite{Na-quiver}で導入された箙多様体に他ならない．

$(\bG,\bN)$の定めるクーロン枝は，$Q$が$ADE$型のとき
に\cite{2016arXiv160403625B}で，アファイン$A$型のとき
に \cite{2016arXiv160602002N}で，別の記述が与えられた．前者の結果を説明
するために，記号の準備を行う．

箙$Q$の向きを忘れ，$ADE$型のディンキン図式とみなし，対応する複素単純群
で随伴型のものを$G$で表す．
\begin{NB}
  ボレル部分群，極大トーラスを取り，$B$, $T$で表す．ワイル群の最長元
  を$w_0$, oppositeボレルを$B_-$とする．
\end{NB}%
$i\in Q_0$に対応する基本余ウェイトを$\varpi_i$, 単純余ルート
を$\alpha_i$と書く．
\footnote{幾何学的佐武対応により，Langlands双対$G^\vee$のウェイト，ルート
  と考えるので，添字$\vee$はあえて省略する．}
上で与えられた $V$, $W$ に対して，$\lambda = \sum (\dim W_i)
\varpi_i$, $\mu = \lambda - \sum (\dim V_i)\alpha_i$とおく．
$G$のアファイン・グラスマン多様体を$\Gr_G$とし，支配的余ウェイ
ト$\lambda$に対応する$G_\cO$軌道を $\Gr_G^\lambda$, その閉包
を$\overline{\Gr}_G^\lambda$で表す．すなわち\subsecref{subsec:geom_satake}
の$\bG$を$G$で置き換えたものである．
さらに $\mu$ も支配的と仮定する．$\mu\le\lambda$ となることか
ら$\Gr_G^\mu\subset\overline{\Gr}_G^\lambda$となる．このと
き\cite{bf14}(\cite{mf}も参照)によ
り，$\Gr_G^\mu$の$\overline{\Gr_G^\lambda}$内の，ある標準的な横断切片が
定義される．すなわち $G[z^{-1}]_1$を，$z=\infty$で値を取る準同
型$G[z^{-1}]\to
G$の核と定義したときに，$G[z^{-1}]_1z^\mu\cap
\overline{\Gr_G^\lambda}$と取る．これを$\oW^\lambda_\mu$で表す．$\mu$が
支配的ではないときの定義は省略するが，そのとき$\oW^\lambda_\mu$は\emph{一般化さ
  れた横断切片}と呼ばれる．
\begin{NB}
余ウェイトを
取り，$\oW^\lambda_{G,\mu}\equiv\oW^\lambda_\mu$を次のデータから成り立
つモジュライ空間として定義する．
\begin{aenume}
\item $\proj^1$上の$G$主束$\scP$
\item $\scP$の$\proj^1\setminus \{0\}$ 上の自明化
  $\sigma\colon \proj^1\setminus \{0\}\times G
  \xrightarrow{\cong} \scP|_{\proj^1\setminus 0}$ であって，$0$での極の次数が
  $\le\lambda$であるもの
\item $\scP$ の $B$簡約$\phi$ であり，次数が $w_0\mu$，$\infty$でのファイバーが
  $\sigma$を通じて$B_-$であるもの
\end{aenume}
このとき(a),(b)は，$\overline{\Gr_G^\lambda}$の点を与える．従って，
写像$\bp\colon\oW^\lambda_\mu\to\overline{\Gr_G^\lambda}$が定義される．
一方，(c)は(b)を通じて$0$自明束の有理的な$B$簡約と思うことができ，点付
き正則写像$f\colon \proj^1\to G/B$
($f(\infty) = B_-$)の全体のなす空間の部分コンパクト化である，zastava空
間 $Z^{-w_0(\lambda-\mu)}$ (\cite{mf})の点を与える．ここ
で，$-w_0(\lambda-\mu)$は写像の次数を表す．したがっ
て，写像 $\bq\colon\oW^\lambda_\mu\to Z^{-w_0(\lambda-\mu)}$が定義される．

特別な場合として，$\lambda=0$のときは，点付き正則写像の全体の空間であり，
また$\mu$が支配的なときには，\cite{bf14}により，$\bp$が局所閉埋め込みに
なり，その像は
$\Gr_G^\mu$の$\overline{\Gr_G^\lambda}$内の(ある標準的な)横断切片になる
ことが知られている．
この理由で，$\oW^\lambda_\mu$ は{\bf 一般化された横断切片}とよばれている．
\end{NB}%

\begin{Theorem}[\protect{\cite{2016arXiv160403625B}}]\label{thm:slice}
  $(\bG,\bN)$ の定めるクーロン枝は $\oW^\lambda_\mu$ と同型である．
\end{Theorem}

クーロン枝の定義自体にアファイン・グラスマンが用いられるが，ここではクー
ロン枝もアファイン・グラスマンと関係しているという結果である．
二つのアファイン・グラスマンは群が異なり，直接的な関係は見られない．


\cite{2017arXiv170900391K}により，\eqref{eq:h:1}のrepelling set
\(
\{ x\in \overline{\Gr}_{G}^\lambda
|
    \lim_{t\to \infty} \nu(t)x = z^\mu
    \}
    \)
は$\oW^\lambda_\mu$のラ
グランジアン部分多様体になる．$\mu$が支配的でない場合
は，$\oW^\lambda_\mu$は$\Gr_G$の部分集合ではないので，正確に
は$\oW^\lambda_\mu$のラグランジアン部分多様体とrepelling setが同型にな
る，という意味である．

\begin{Remark}\label{rem:symmetrizable}
  上のクーロン枝の定義から現れる$G$は，$ADE$型，あるいはより一般の
  箙を考えたとしても，対称なKac-Moodyリー環に対応する群にな
  る．\cite{2019arXiv190706552N}では，クーロン枝の定義を修正することに
  より，$BCFG$型や，より一般の対称化可能なKac-Moodyリー環に対応するもの
  を実現した．ポイントは $\bG$ のアファイン・グラスマン $\Gr_\bG$ は，
  一般線形群のアファイン・グラスマンの積 $\prod_i \Gr_{\GL(V_i)}$ であ
  るが，各頂点 $i$ ごとに異なる形式的円盤を考えることである．すなわ
  ち $i$ごとに変数$z_i$ を用意して，形式的円盤 $D_i = \Spec
  \CC[[z_i]]$ を考える．カルタン行列 $C=(a_{ij})$ が対称化可能であると
  は，各$i$に対して正整数 $d_i$ が存在して $d_ia_{ij} = d_j a_{ji}$ と
  なるときをいうが，これに対応して $z_i$ の $d_i$ 次被覆
  を $z_i^{1/d_i} = z_j^{1/d_j}$ のように関係させる．そこ
  で $\Hom(V_i,V_j)$ に対応する $\cR$ の成分を，被覆のアファイン・グラ
  スマン上で考えることにより，クーロン枝の定義を変更する．定理~\ref{thm:slice}
  は，$G$が$BCFG$型でも成立する．
\end{Remark}

今の場合は量子化されたクーロン枝$\cAh$は，シフトされたヤンギアンと呼ば
れる代数の商になる．(\cite{2016arXiv160403625B}のAppendixを参照．)
詳細は省略するが，定理~\ref{thm:classical}のあとに述べた$\cAh$ の有理関
数係数差分作用素による実現を用いて，$\cAh$の生成元がシフトされたヤンギ
アンの定義関係式を満たすことにより証明される．

\begin{NB}
\subsection{シフトされたヤンギアン}

次に，\cite{2016arXiv160403625B}の付録で決定された，$ADE$型箙ゲージ理論
の量子化されたクーロン枝について，主張を述べるために記号を準備する．

$\fg$を$Q$の定める複素単純リー環とする．上のように$\mu$を$\fg$の余ウェイトとする．
\begin{Definition}
  \textup{(1)}
  $Y_\infty(\fg)$ は，生成元 $ E_i^{(q)}, F_i^{(q)}, H_i^{(p)} $
  ($ q \in \ZZ_{> 0}$, $ p \in \ZZ $, $ i \in Q_0 $) と関係式
  \begin{align*}
    [H_i^{(p)}, H_j^{(p')}] &= 0,  \\
    [E_i^{(p)}, F_j^{(q)}] &=   \delta_{ij} H_i^{(p+q-1)}, \\
    [H_i^{(p+1)},E_j^{(q)}] - [H_i^{(p)}, E_j^{(q+1)}] &= \frac{ \alpha_i \cdot \alpha_j}{2} (H_i^{(p)} E_j^{(q)} + E_j^{(q)} H_i^{(p)}) , \\
    [H_i^{(p+1)},F_j^{(q)}] - [H_i^{(p)}, F_j^{(q+1)}] &= -\frac{ \alpha_i \cdot \alpha_j}{2} (H_i^{(p)} F_j^{(q)} + F_j^{(q)} H_i^{(p)}) , \\
    [E_i^{(p+1)}, E_j^{(q)}] - [E_i^{(p)}, E_j^{(q+1)}] &= \frac{ \alpha_i \cdot \alpha_j}{2} (E_i^{(p)} E_j^{(q)} + E_j^{(q)} E_i^{(p)}), \\
    [F_i^{(p+1)}, F_j^{(q)}] - [F_i^{(p)}, F_j^{(q+1)}] &= -\frac{ \alpha_i \cdot \alpha_j}{2} (F_i^{(p)} F_j^{(q)} + F_j^{(p)} F_i^{(q)}),\\
    i \neq j, N = 1 - \alpha_i \cdot \alpha_j \Rightarrow
    \operatorname{sym} &[E_i^{(p_1)}, [E_i^{(p_2)}, \cdots [E_i^{(p_N)}, E_j^{(q)}]\cdots]] = 0, \\
    i \neq j, N = 1 - \alpha_i \cdot \alpha_j \Rightarrow
    \operatorname{sym} &[F_i^{(p_1)}, [F_i^{(p_2)}, \cdots
                         [F_i^{(p_N)}, F_j^{(q)}]\cdots]] = 0
  \end{align*}
  で定義された$\CC$代数である．ここ
  で，$\operatorname{sym}$は，$p_1,\dots,p_N$に関する対称化を表す．

  \textup{(2)} シフトされたヤンギア
  ン $Y_\mu(\fg)$ を，$Y_\infty(\fg)$にさらに条件
  \begin{equation*}
    H_i^{(p)} = 0 \quad p < -\langle\mu,\alpha_i\rangle,\qquad
    H_i^{(\langle\mu,\alpha_i\rangle)} = 1
  \end{equation*}
  を課して定められる商代数として定義する．
\end{Definition}

$\mu = 0$ のときは，インデックスが一つづつずらされていることを除
き，Drinfeldのヤンギアンに'new presentation' \cite{MR914215} のもとで等
しい．また，$\mu$が支配的なときには，$H_i^{(p)}$と$F_i^{(q)}$ のインデッ
クスを$\langle\mu,\alpha_i\rangle$だけずらすことによ
り，$Y_\mu(\fg)$は$Y_0(\fg)$の中の部分代数として埋め込め，これ
は\cite{kwy}で定義されたシフトされたヤンギアンと同じになる．

\begin{Theorem}[\protect{\cite[App.~B]{2016arXiv160403625B}}]
  (具体的に記述される)環準同型 $Y_\mu(\fg)\to
  \left.\cAh\right|_{\hbar=1}$が存在し，$\mu$が支配的なときは全射である．
\end{Theorem}

なお，$\mu$が支配的である，という条件は \cite{2019arXiv190307734W} で外された．

証明には，\secref{sec:local}の局所化定理の，量子化されたクーロン枝版を用
いる．すなわち，$\cAh$は，$\mathfrak t$上の有理差分作用素のなす環に埋め
込むことができる．一方，$Y_\mu(\fg)$については，有理差分作用素を用いて
表現を実現できることが，$\mu=0$のときは\cite{GKLO}により，$\mu$が支配的
なときは\cite{kwy}によって示されており，一般の$\mu$に拡張することもでき
る．そこで，シフトされたヤンギアンの生成元に対応する差分作用素が，量子
化されたクーロン枝に属していることをみれば，環準同型の存在が従う．全射
であることは，$\hbar=0$に誘導される準同型が全射であることを示せばよ
く，$\mu$が支配的なときは，左辺は横断切片の量子化であること
が\cite{kwy}で示されているので，これが従う．
\end{NB}%

\subsection{ジョルダン箙とDAHA} 

$Q$をジョルダン箙とする．すなわち，頂点は一つで，その頂点を自分自身に結ぶ
辺が一つあるものである．対応するゲージ理論は，$V$, $W$を有限次元複素ベクトル空間として
$\bG = \GL(V)$, $\bN = \End(V)\oplus\Hom(W,V)$と取ったものである．

\begin{Theorem}[\protect{\cite[Prop.~3.24]{2016arXiv160403625B}}]
$\bG = \GL(V)$, $\bN = \End(V)\oplus\Hom(W,V)$に対応するクーロン枝は，
曲面$\mathcal S_\ell = \{ (x,y,z)\in\CC^3 \mid xy = z^\ell\}$ の$n$次対称積
$S^n(\mathcal S_\ell)$に等しい．ただし，$\ell = \dim W$, $n=\dim V$である．
\end{Theorem}

$W = 0$ ($\ell = 0$)
のときは，$\bN = \mathfrak{gl}(V)$なので，\subsecref{subsec:quantum}の
最後に出てきた例の，$\bG =
\GL(V)$の場合であり，$\mathcal S_0 = \CC\times\CC^\times$なの
で，$\mathfrak t\times \bT^\vee/\mathbf W$と同じになる．そこで言及した
ように，対応する量子化されたクーロン枝は，三角DAHAのsphereical partになる．
一般の$W$については，その拡張として次が分かる．

\begin{Theorem}[\cite{2016arXiv160800875K}]
  $\ell \ge 1$ のとき，上の$\bG$, $\bN$に対する量子化されたクーロン枝は，
  wreath積 $S_n\ltimes (\ZZ/\ell\ZZ)^n\subset\grpSp(n)$ に付随した
        有理Cherednik代数のspherical part
        に
        同型である．
\end{Theorem}

証明は，シフトされたヤンギアンの場合と同様に，有理差分作用素の環への埋
め込みを用いる．DAHAについての対応する埋め込みは，Demazure-Lusztig作用
素で実現される．

\section{量子化されたクーロン枝の表現}

\subsection{c-ウェイト加群}

量子化されたクーロン枝$\cAh$ は，可換な部分代
数 $H^*_{\bG\times\CC^\times}(\mathrm{pt})$を含んでいたことを思い出そう．
以下，$\hbar = 1$と特殊化した $\cAh|_{\hbar=1}$ を考える．これを，$\cA$で表す．

$H^*_{\bG}(\mathrm{pt}) = \CC[\ft]^{\mathbf W}$に注意しよう．

\begin{Definition}
  $\cA$の表現 $M$を考える．

  (1) $M$の{\bf
    c-ウェイト空間}とは，$H^*_\bG(\mathrm{pt})$に関する同時一般化固有空間
  であって，$0$でないものをいう．このとき，同時固有
  値
  $\lambda\in\Hom_{\mathrm{alg}}(H^*_\bG(\mathrm{pt}) , \CC)
  =\mathrm{Specm}(H^*_\bG(\mathrm{pt})) = \ft/{\mathbf W}$
  を$\ft$にリフトしたものを{\bf c-ウェイト}という．

  (2) $M$が{\bf c-ウェイト加群}であるとは，$M$はc-ウェイト空間の直和に
  分解するときをいう．
\end{Definition}

同時固有空間になっていることまでは課していないので，一般
に$H^*_{\bG}(\mathrm{pt})$の元 $f$ は，c-ウェイト空間に
$f(\lambda) + \text{nilpotent}$で作用していることに注意する．

また，$\bG$の余ウェイトと$M$のc-ウェイトを，共にウェイトをつけて呼ぶの
は混乱を招きかねないように思うが，他の用語を思いつかないので，とりあえ
ず後者には c をつけて区別していることで十分であると諦めることにする．

\begin{Lemma}
    $M$が直既約であれば，$M$がc-ウェイト $\lambda$, $\mu$ を持つとする
    と，$\mu - \lambda$, 正確にはそれぞれを$\ft$にリフトしたものの差
    は，$\bT$の余ウェイトである．
\end{Lemma}

以下では，簡単のために $\lambda$ 自身が$\bT$の支配的余ウェイトである場合
を考える．このとき他の
c-ウェイトは，やはり$\bT$の支配的余ウェイトである．この条件を外しても，同
様の考察は可能であるが，次の節の固定点集合の記述が面倒になる．

$\lambda$が与えられたとき，$\tilde\lambda\colon\CC^\times\to
\bT\times\CC^\times$を$\tilde\lambda(\tau)
= (\lambda(\tau),\tau)$ によって定義する．

\subsection{固定点集合}

以下では，$\cA$のc-ウェイト加群 $M$ を幾何学的に解析する理論を紹介する．
これは，Ginzburgの合成積代数の解析の理論 \cite{CG}を，無限次元多様体を
用いるDAHAの場合に修正した\cite{MR3013034}の枠組みに基づくものである．

この理論に現れるのは，
$\cT = \bG_\cK\times^{\bG_\cO}\bN_\cO$, $\bN_\cK$ の
$\tilde\lambda(\CC^\times)$固定点集合
$\cT^{\tilde\lambda}$, $\bN_\cK^{\tilde\lambda}$ と，
射$\Pi\colon\cT\to\bN_\cK$の制限
$\Pi^{\tilde\lambda}\colon 
\cT^{\tilde\lambda} \to \bN_\cK^{\tilde\lambda}$ である．

\begin{NB}
$\Gr_\bG$ への$\bT\times\CC^\times$作用を考える．ここで，$\bT$ は
$\bG_\cO$ の部分群として作用する．
$\bT$の余指標$\lambda$と$m\in\ZZ_{>0}$に対して，$\bT\times\CC^\times$へ
の群準同型$\CC^\times\to \bT\times\CC^\times$
$\tau\mapsto (\lambda(\tau),\tau^m)$を考える．
以下では，簡単のために$m=1$の場合のみを考える．
\end{NB}%

\begin{Lemma}
    $\tilde\lambda$ に関する $\Gr_\bG$の固定点集合は，
    \begin{equation*}
        \bigsqcup_{\mu} \{\lambda(z)^{-1} g [z^\mu] \in\Gr_{\bG}
        \mid {g\in\bG} \}
    \end{equation*}
    となる．ここで，$\mu$は$\bG$の支配的な余ウェイトであ
    る．$z^{\mu-\lambda}$を含む連結成分は，$\mu$に対応する放物部分群
    を$\mathbf P_\mu$として，部分旗多様体 $\bG/\mathbf P_\mu$に同型である．
\end{Lemma}

\begin{proof}
    $\Gr_\bG\cong\Omega\bG_c$を用いると，$c(z)\in\Omega\bG_c$が固定され
    るのは，
    \begin{equation*}
        \lambda(\tau) c(z\tau) c(\tau)^{-1} \lambda(\tau)^{-1}
        = c(z)
    \end{equation*}
    が，$z\in S^1$, $\tau\in S^1$について成り立つときである．これ
    は，$z\mapsto \lambda(z)c(z)$ が
    $S^1\to \bG_c$として群準同型であることに他ならない．群準同型は，支
    配的な $\mu\colon S^1\to\bT_c$と
    $\bG$共役になるので，$c(z) = \lambda(z)^{-1} g \mu(z) g^{-1}$ とな
    る．また $\mu$は，uniqueである．これをアファイン・グラスマンの記法
    に直すと，上の主張になる．後半は，$g[z^\mu] = [z^\mu]$ となる
    $g$の全体を求めればよく，これは$\mathbf P_\mu$に他ならない．
\end{proof}

\begin{NB}
    $g[z^\mu] = [z^\mu]$となる条件を書いてみると，
    $\mathbf P_\mu = \{ g\in\bG \mid z^{-\mu} g z^\mu \in \bG_\cO\}$である．
\end{NB}%

次に，$\cT$における固定点集合を求める．
まず，$\bN_\cK$における固定点を求めてみる．作用は
$s(z)\mapsto \lambda(\tau) s(z\tau)$ で与えられるから，
\begin{equation*}
    \bN_\cK^{\tilde\lambda}
    \cong
    \bN
\end{equation*}
である．ここで，右から左への同型は $\bN\ni s \mapsto
s(z) = z^{-\lambda} s$ で与えられる．

上の補題において$g=1$に対応する$\Gr_\bG$の固定点上のファイバーを考える
と，
\begin{NB}
\[
    [z^{\mu-\lambda},s(z)] \mapsto
    [\lambda(\tau) (z\tau)^{\mu-\lambda}, s(z\tau)]
    = [z^{\mu-\lambda}, \tau^\mu s(z\tau)]
\]
\end{NB}%
作用が$\lambda(\tau)s(z\tau)$の代わりに $\mu(\tau)s(z\tau)$となることに
注意して，
\begin{equation*}
    \{ s(z)\in\bN_\cO \mid s(z\tau) = \tau^{-\mu} s(z) \}
    \cong
    \{ s\in\bN \mid z^{-\mu} s \in \bN_\cO \}
\end{equation*}
である．ここで，右から左への同型は，上と同様に $s(z) = z^{-\mu} s$ で与
えられる．
\begin{NB}
    $s(z)$は $z=1$のときの値，$s\in\bN$ で決まる．
\end{NB}%
この$\bN$の部分空間を，$\bN^\mu_{\le 0}$
とおく．$\bN$を$\mu$に関してウェイト空間分解したときに，ウェイトが非正
のものの直和に他ならない．これは，$\mathbf P_\mu$不変である．
\begin{NB}
    実際，$g\in \mathbf P_\mu$
    に対して，$z^{-\mu} g s = z^{-\mu} g z^{\mu} z^{-\mu} s$であ
    り，$z^{-\mu} g z^\mu\in\bG_\cO$, $z^{-\mu}s\in\bN_\cO$に注意して，
    $z^{-\mu} g s\in\bN_\cO$である．
\end{NB}%
従って
\begin{Lemma}
    $\tilde\lambda$ に関する固定点集合は，

    (1) $\bN_\cK$については，$\bN$と同型である．

    (2)
    $\cT$については，$\bigsqcup_\mu \bG\times^{\mathbf P_\mu}
    \bN^\mu_{\le 0}$
    と同型になる．$\Pi\colon\cT\to\bN_\cK$の制限は，自然な写
    像 $\bG\times^{\mathbf P_\mu}\bN^\mu_{\le 0}\to \bN$;
    $[g,s]\mapsto gs$である．
\end{Lemma}

\begin{NB}
    (2) の最後の主張をチェックする．

    $z^{-\lambda} g [z^\mu]$
    におけるファイバーは，$[z^{-\lambda} g z^\mu, z^{-\mu}s]$
    ($s\in\bN^{\mu}_{\le 0}$) である．これを$\Pi$で送る
    と，$z^{-\lambda}g s\in \bN^{\tilde\lambda}_\cK$ であり，上の同型
    で$\bN$ の元と思ったときは，$gs$ である．
\end{NB}

$\bG\times^{\mathbf P_\mu} \bN^\mu_{\le 0}$ を $\cT^\mu$ とおく．
この設定で，合成積代数
\begin{equation}\label{eq:6}
    \prod_\mu \bigoplus_{\mu'}
    H_*(\cT^\mu\times_{\bN}\cT^{\mu'})
\end{equation}
を考えることができる．($\mu$, $\mu'$について無限和なので，単位元を持
つように，$\prod$, $\oplus$ を取る必要がある．)

\begin{Theorem}[証明の詳細は\cite{modules}で発表予定]
    (1) 環準同型 $\cA\to 
        \prod_\mu \bigoplus_{\mu'}
    H_*(\cT^\mu\times_{\bN}\cT^{\mu'})$ が，幾何学的に構成できる．

    (2) c-ウェイト $\lambda$ を持つ直既約 c-ウェイト加群 $M$ は，
    上の環準同型を通じて得られる．
\end{Theorem}

従って，c-ウェイト $\lambda$ を持つ c-ウェイト加群の表現論を調べる
には，\eqref{eq:6}の表現論を調べればよいことになる．

\subsecref{subsec:triples}で紹介したよう
に，$\cT^\mu\times_{\bN}\cT^{\mu'}$は Steinberg多様体の類似物であり，合
成積代数はGinzburgの理論 \cite{CG}を用いて解析することができる．たとえば，
次が分かる．
\begin{Theorem}
    c-ウェイト $\lambda$ を持つ既約 c-ウェイト加群 $M$ は，
    射 $\cT^\mu\to\bN$
    による$\cT^\mu$の定数層のpushforwardの分解の中にシフトを除いて現れ
    る単純偏屈層の同型類と一対一に対応する．ただし，後者では$\bG$の支配
    的余ウェイト$\mu$をすべて動かす．
\end{Theorem}

単純偏屈層をすべて決定することは，一般には難しい．しかし，特別なときは，
すでに調べられている空間が固定点集合として出てくることがありえる．たと
えば，\secref{sec:quiver}で$W=0$の場合には，固定点集合は本質的にLusztigが
標準基底の定義に用いた空間\cite{Lu-can2}に他ならない．従って，合成
積代数は\cite{MR2837011}より，KLR代数である．よって既約表現は，双対標準
基底の元と一対一に対応する．

ジョルダン箙の場合には，$\lambda$として余ウェイトになっているものだけで
はなく，より一般に，$m\in\ZZ_{>0}$に対して $m\lambda$ が余ウェイトになっ
ているものを考えるべきである．この場合は，$\tilde\lambda(\tau) =
(m\lambda(\tau),\tau^m)$に関する固定点集合を考える必要があり，す
ると頂点が$m$個の巡回箙がジョルダン箙に代わって現れる．この場合はアファイ
ン $A_{m-1}$型の双対標準基底との対応が従う．

\begin{NB}
\section{Geometric Satake correspondence for complex reductive groups\\\hfill --- disclaimer}

Our main result \ref{thm:satake} can be regarded as an affine Lie
algebra $\algsl(n)_{\mathrm{aff}}$ version of the geometric Satake
correspondence for a complex reductive group $\bG$
\cite{Lus-ast,1995alg.geom.11007G,Beilinson-Drinfeld,MV2}. Our
approach is close to one due to Mirkovi\'c-Vilonen
\cite{MV2}. Fortunately we have many good expository articles on
\cite{MV2}, and there is no reason to try to produce a worse one.
There is another reason to decide not to explain geometric Satake
correspondence: We avoid to use sheaf theoretic language as much as
possible in the spirit of the talk in Algebraic Geometry Symposium.
If we would explain geometric Satake, this decision makes no sense.

We also omit historical accounts. Interested readers should read
\cite[\S3(viii)]{2016arXiv160403625B} and \cite{fnkl_icm}.
\end{NB}%

\section{箙多様体とKac-Moodyリー環}

\secref{sec:quiver}の箙ゲージ理論の設定において，\subsecref{subsec:Higgs}のヒッ
グス枝を考えたものが，箙多様体にほかならない．箙多様体による
Kac-Moodyリー環の表現の構成については，以前に\cite{MR1802956}で
解説したが，幾何学的佐武対応との対比を明確にするために復習する．

\subsection{定義}

$\bG$, $\bN$ を \eqref{eq:7}
のように取る．$\bM = \bN\oplus\bN^*$とする．\subsecref{subsec:Higgs}で注意し
たように，圏論的商を用いたハミルトニアン簡約 $\mu^{-1}(0)\dslash
\bG$ と，$\bG$の乗法的指標 $\chi$ を取り，幾何学的不変式論によって取っ
た商を用いる $\mu^{-1}(0)\dslash_\chi \bG$ が定義される．自然な射影的
射$\pi\colon \mu^{-1}(0)\dslash_\chi \bG \to \mu^{-1}(0)\dslash \bG$ が
定まる．
$\chi$ が generic であれば $\mu^{-1}(0)\dslash_\chi \bG$ は滑らかであり，
また多くの場合に $\pi$ は特異点解消になることが知られている．


\subsection{箙多様体と可積分最高ウェイト表現}\label{subsec:quiver_integrable}

さらに 箙 $Q$ は，頂点を自分自身に結ぶループを持たないと仮定す
る．\secref{sec:quiver}と同様に箙の辺の向きを忘れてディンキン図式を考え，
対応する対称なKac-Moodyリー環を $\mathfrak g$ で表す．
$Q_0$-graded なベクトル空間 $V$, $W$ に対して余ウェイ
ト $\lambda$, $\mu$ を定理~\ref{thm:slice}のように定める．以下では，しばら
く$\mathfrak g$ は対称であると仮定するので，Langlands双対
も $\mathfrak g$ となり，$\lambda$, $\mu$ は $\mathfrak g$のウェイトで
あるとも考える．
さらに二通りのハミルトニアン簡約 $\mu^{-1}(0)\dslash\bG$,
$\mu^{-1}(0)\dslash_\chi\bG$ をそれぞれ $\fM_0(\lambda,\mu)$,
$\fM_\chi(\lambda,\mu)$ で表す．さらに $\chi$ は $\bG = \prod
\GL(V_i)$ の行列式の積で与えられるとする．

\begin{Theorem}[\cite{Na-quiver,Na-alg}]\label{thm:quiver}
  \textup{(1)} $0$ を $0\in \bN\oplus\bN^*$ に対応す
  る $\fM_0(\lambda,\mu)$の点とする．このと
  き $\fL_\chi(\lambda,\mu)\defeq
  \pi^{-1}(0)$ は $\fM_\chi(\lambda,\mu)$ の(一般には特異点を持つ)ラグ
  ランジアン部分多様体である．

  \textup{(2)} $\fL_\chi(\lambda,\mu)$ の最高次のホモロジーの直和
  \begin{equation*}
    \bigoplus_\mu H_{\operatorname{top}}(\fL_\chi(\lambda,\mu))
  \end{equation*}
  は，$\mathfrak g$ の可積分最高ウェイト表現の構造を持つ．
\end{Theorem}

$\mathfrak g$の可積分最高ウェイト表現は自動的に既約であり，整な支配的ウェ
イトが最高ウェイトとなってパラメトライズされる．$\lambda$ に対応する可
積分最高ウェイト表現を $V(\lambda)$ で表す．ウェイト分解を持ち，
$V_\mu(\lambda)$ をウェイト $\mu$ のウェイト空間とする．上の結果に
おい
て$H_{\operatorname{top}}(\fL_\chi(\lambda,\mu))$ は $V_\mu(\lambda)$
に対応する．

$H_{\operatorname{top}}(\fL_\chi(\lambda,\mu))$ は，
$\fL_\chi(\lambda,\mu)$ の既約成分の基本類が与える基底を持つ．また，既
約成分の集合の和 $\bigsqcup_\mu
\operatorname{Irr}\fL_\chi(\lambda,\mu)$ に，柏原の意味のクリスタルの構
造を入れることができ，さらに量子展開環の表現のクリスタルと同型になるこ
とが知られている\cite{KS,Saito,Na-Tensor}．
ここで，柏原の意味のクリスタルの定義や，この結果の正確な主張は説明しな
い．雑にいうと，$\fL_\chi(\lambda,\mu)$ のすべての既約成分を最高ウェイ
トベクトル $\fL_\chi(\lambda,0) = \{0\}$ から出発して，$\mu$ を動かしな
がら再帰的に構成することができることを主張している．また，この構造が上
の定理において，表現が最高ウェイトであることの証明に使われる．

あとで解説する予想~\ref{thm:satake} と，上の定理~\ref{thm:quiver}とにお
ける，$\mathfrak g$の表現の構成の類似性を明確にするために，表現がどのよ
うに定義されるのかを解説しよう．

Kac-Moodyリー環 $\mathfrak g$ はディンキン図式の頂点 $i$ に
対応した生成元 $e_i$, $f_i$ および, 可換なカルタン部分環 $\mathfrak h$
で生成され，ある関係式で定義される．
$\mathfrak h$の表現は，ウェイト $\mu$，すなわ
ち $H_{\mathrm{top}}(\fL_\chi(\lambda,\mu))$ の $V_i$, $W_i$ の次元で決
まる．

一方， $e_i$, $f_i$ は
$\fM_\chi(\lambda,\mu)\times
\fM_\chi(\lambda,\mu-\alpha_i)$ 内のcorrespondence で与えられる．ここで，
ウェイト $\mu-\alpha_i$ は，$Q_0$-graded ベクトル空間 $V\oplus S_i$から
来ることに注意する．ただし $S_i$ は頂点 $i\in Q_0$ にのみ$1$次元ベクト
ル空間を持ち，他の頂点では $0$ となる $Q_0$-graded ベクトル空間である．

さらに $\chi_i\colon G\to \CC^\times$ を$j\neq i$ の $\GL(V_j)$ の行
列式の積として与えられる乗法的指標とする．
これに対応する幾何学的不変式論的な商 $\fM_{\chi_i}(\lambda,\mu)$,
$\fM_{\chi_i}(\lambda,\mu-\alpha_i)$ を取る．射影的
射 $\pi$
は，$\fM_\chi(\lambda,\mu)\to \fM_{\chi_i}(\lambda,\mu)\to
\fM_0(\lambda,\mu)$ のように分解
し，$\fM_{\chi_i}(\lambda,\mu-\alpha_i)$ についても同様である．
さらに
閉埋め込み$\fM_{\chi_i}(\lambda,\mu)
\to\fM_{\chi_i}(\lambda,\mu-\alpha_i)$ が存在する．これは，$V,W$に対し
て $\mu^{-1}(0)$ に入る $\bN\oplus\bN^*$ の元を，$S_i$成分を $0$にして
拡張して，$V\oplus S_i,W$ に対応する元と見ることにより定められる．
指標 $\chi$ の幾何学的不変式論の安定性の条件は，この拡張の操作で保たれ
ないので，$\fM_\chi(\lambda,\mu)\to\fM_\chi(\lambda,\mu-\alpha_i)$とい
う閉埋め込みは定義されないが，$\chi_i$ に置き換えると，well-defined に
なる．これは，$\chi_i$ の定め方により，$\GL(V_i)$ については圏論的な商
になっていることによる．

そこで，ファイバー積
\begin{equation*}
  \fM_\chi(\lambda,\mu)\times_{\fM_{\chi_i}(\lambda,\mu-\alpha_i)}
  \fM_\chi(\lambda,\mu-\alpha_i),
\end{equation*}
を考える．ただし，$\fM_\chi(\lambda,\mu)\to \fM_{\chi_i}(\lambda,\mu-\alpha_i)$
は，$\fM_\chi(\lambda,\mu)\to
\fM_{\chi_i}(\lambda,\mu)\to\fM_{\chi_i}(\lambda,\mu-\alpha_i)$の合成である．
これが$\fM_\chi(\lambda,\mu)\times \fM_\chi(\lambda,\mu-\alpha_i)$内の
ラグランジアン部分多様体であることが知られている．

このファイバー積の中に，一つの既約成分 $\mathfrak P_i(\lambda,\mu)$ が
ある．それは，おおざっぱにいうと $V\oplus S_i$ に対応する
元 $x'$ を $V$ に制限してできる $V$ に対応してできる元を $x$ として，
組 $(x,x')$ で与えられるものである．ただし，$\GL(V_i\oplus S_i)$ の群の
作用まで込めて考えないといけないので，$V$ への制限をどのように定義する
のかは慎重に取り扱う必要がある．その技術的な詳細は省略する．

さらに，各成分への射影
$p_1, p_2\colon \mathfrak P_i(\lambda,\mu)\to \fM_\chi(\lambda,\mu)$,
$\fM_\chi(\lambda,\mu-\alpha_i)$ は固有であり，従って, 線形写像
\begin{equation*}
  H_{\mathrm{top}}(\fL_\chi(\lambda,\mu)) \leftrightarrows
  H_{\mathrm{top}}(\fL_\chi(\lambda,\mu-\alpha_i))
\end{equation*}
が，合成積 $p_{2*}p_1^*$, $p_{1*} p_2^*$ によって定義される．符号を除い
て，これが生成元の $e_i$, $f_i$ の定義である．

\subsection{箙多様体によるテンソル積}\label{subsec:tensor1}

定理~\ref{thm:quiver} において構成されるのは，既約表現である．次
に，Lusztig \cite{MR1714628} と Varagnolo-Vasserot \cite{VV-std} の先行
研究に動機付けされて考えた，テンソル積表現の構成 \cite{Na-Tensor} につ
いて解説する．

定理~\ref{thm:quiver}において最高ウェイト $\lambda$ は $(\dim
W_i)_{i\in Q_0}$ により与えられていたことを思い出そ
う．$Q_0$-graded ベクトル空間の直和分解 $W = W^1\oplus W^2$ を考え
る．$(\dim W^1_i)$ と $(\dim W^2_i)$により，二つの整な支配的ウェイ
ト $\lambda^1$, $\lambda^2$ が対応する．以下では
$V(\lambda^1)\otimes V(\lambda^2)$ を構成する．構成自体は，$W^1$,
$W^2$ を入れ替えると一見違うものになるが，テンソル
積 $V(\lambda^1)\otimes V(\lambda^2)$ は成分を入れ替えても同型になるこ
とに注意しよう．表現を量子ループ代数 $\mathbf U_q(\mathbf L\mathfrak
g)$ に持ち上げると，テンソル積の成分の順序に依存するので，幾何学的な構
成が入れ替えについて対称ではないのは自然である．

1パラメータ変換群 $\nu\colon\CC^\times\to
\GL(W)$を$\nu(t) = \operatorname{id}_{W^1}\oplus
t\operatorname{id}_{W^2}$ によって定義する．
これを通じて $\CC^\times$ が $\fM_0(\lambda,\mu)$,
$\fM_\chi(\lambda,\mu)$ に自然に作用する．
そこで $\fM_0(\lambda,\mu)$ の attracting set
\begin{equation}\label{eq:91}
\begin{gathered}
  \fT^{\nu}_0(\lambda,\mu) \defeq \left\{ x\in\fM_0(\lambda,\mu)
    \,\middle|\,
     \text{$\lim_{t\to 0} \nu(t)x$ が存在する}
  \right\}, \\
  \widetilde{\fT}^{\nu}_0(\lambda,\mu) \defeq \left\{ x\in\fM_0(\lambda,\mu) \,\middle|\,
    \lim_{t\to 0} \nu(t)x = 0\right\}.
\end{gathered}
\end{equation}
を考え，その$\pi$ による逆像を，それぞ
れ $\fT^{\nu}_\chi(\lambda,\mu)$,
$\widetilde{\fT}^{\nu}_\chi(\lambda,\mu)$ で表す．

\begin{Theorem}[\cite{Na-Tensor}]\label{thm:tensor}
  \textup{(1)} $\widetilde{\fT}^{\nu}_\chi(\lambda,\mu)$ は
   $\fM_\chi(\lambda,\mu)$ のラグランジアン部分多様体である．

   \textup{(2)} $\widetilde{\fT}^{\nu}_\chi(\lambda,\mu)$ の最
   高次のホモロジーの直和
  \begin{equation*}
    \bigoplus_\mu H_{\mathrm{top}}(\widetilde{\fT}^{\nu}_\chi(\lambda,\mu))
  \end{equation*}
  は $\mathfrak g$ の可積分表現の構造を持ち，二つの最高ウェイト表現のテンソル積
  $V(\lambda^1)\otimes V(\lambda^2)$ と同型である．
\end{Theorem}

$\mathfrak g$ の表現の構成は，定理~\ref{thm:quiver}と同様に合成積で与えられ
る．表現が正しい`大きさ' を持つことを確認するため
に， $\widetilde{\fT}_\chi^{\nu}(\lambda,\mu)$ の既約成分を記述
しよう．
まず，$\nu$固定点集合 $\fM_\chi(\lambda,\mu)^{\nu}$ は
\begin{equation*}
  \fM_\chi(\lambda,\mu)^{\nu} \cong \bigsqcup_{\mu=\mu^1+ \mu^2}
  \fM_\chi(\lambda^1,\mu^1)\times \fM_\chi(\lambda^2,\mu^2)
\end{equation*}
のように分解することを注意しよう．右辺から左辺への射は，`直和'を取るこ
とで証明される．さらに，同型写像であることは，$\fM_\chi(\lambda,\mu)$が
滑らかで，fine なモジュライ空間であることを使って証明される．詳しく
は \cite[Lemma~3.2]{Na-Tensor} を参照．
$\widetilde{\fT}^{\nu}_\chi(\lambda,\mu)$ の点について，極
限 $\lim_{t\to 0} \nu(t)x$ が存在し，それは $\nu$固定点である．上の
分解に従って，次の分解が誘導される．
\begin{equation*}
  \widetilde{\fT}^{\nu}_\chi(\lambda,\mu) = \bigsqcup_{\mu=\mu^1+ \mu^2}
  \widetilde{\fT}^{\nu}_\chi(\lambda^1,\mu^1;\lambda^2,\mu^2).
\end{equation*}
すると $\widetilde{\fT}^{\nu}_\chi(\lambda^1,\mu^1;\lambda^2,\mu^2)$ は
$\fL_\chi(\lambda^1,\mu^1)\times\fL_\chi(\lambda^2,\mu^2)$ 上のベクトル束であり，
底空間への射影は，極限 $\lim_{t\to 0} \nu(t)x$ で与えられる．
これは，元 $x$ が短完全列$0\to x^2 \to x \to x^1\to
0$ で $x^1\in\fL_\chi(\lambda^1,\mu^1)$,
$x^2\in\fL_\chi(\lambda^2,\mu^2)$ となるものに対応し，ファイバーのベク
トル空間の構造は，$\operatorname{Ext}^1$ から来
る．\cite[Remark~3.16]{Na-Tensor}を参照．
このとき，$\widetilde{\fT}^{\nu}_\chi(\lambda,\mu)$ の既約成分は
$\widetilde{\fT}^{\nu}_\chi(\lambda^1,\mu^1;\lambda^2,\mu^2)$の閉包であり，
従って, 既約成分は
\(
\bigsqcup_{\mu=\mu^1+\mu^2} \operatorname{Irr}\fL_\chi(\lambda^1,\mu^1)\times
\operatorname{Irr}\fL_\chi(\lambda^2,\mu^2)
\)
でパラメトライズされ，確かにテンソル積と同じ大きさを持っている．

例を与えよう．$Q$として $A_{n-1}$ 型のものを取り，$V = \bigoplus V_i$ を $V_i = \CC$
($1\le i\le n-1$), $W = \bigoplus W_i$ を $W_i = \CC$
($i = 1, n-1$), $W_i = 0$ ($1 < i < n-1$) と取る．($n=2$のときは
$W_1 = \CC^2$と取る．) 対応する箙多様体
$\fM_0(\lambda,\mu)$ は $\CC^2/(\ZZ/n\ZZ)$, すなわち $A_{n-1}$型の単純特異点であり， $\fM_\chi(\lambda,\mu)$ はその極小特異点解消である．ラグランジアン部分多様体 $\fL_\chi(\lambda,\mu)$ は $n-1$個の複素射影直線のチェーンである．そこで $W = W^1\oplus W^2$ 
($\dim W^1=\dim W^2 = 1$)という分解を取る．
固定点集合 $\fM_\chi(\lambda,\mu)^{\nu}$
は，$n$個の孤立した点であり，射影直線の北極と南極であり，そのうち
の $(n-2)$点は射影直線が交叉しているところで，残りの二点はチェイン
の両端にある．
$\widetilde\fT^{\nu}_\chi(\lambda,\mu)$ は，もう一つ余分な既約
成分を持つ．分解
$W=W^1\oplus W^2$の取り方に応じて，両端のいずれかの点を通る直線が付け
加えられる．図~\ref{fig:tensor}では，点線で表されているのが付け加わっ
た既約成分である．

対応する表現は，$\algsl_n$のベクトル表現 $\CC^n$ と双対表
現 $(\CC^n)^*$ のテンソル積 $\CC^n \otimes (\CC^n)^*$ である．それ
は $\algsl_n$ の随伴表現と自明表現の直和に分解するが，付け加えられた既
約成分が自明表現に対応する．

\begin{figure}[htbp]
    \centering
\begin{tikzpicture}[scale=0.3]
\draw (0,10) parabola bend (15,-5) (30,10);
\draw[thick,fill] (9,0) circle (0.3);
\draw[thick,bend left,distance=40] (9,0) to (12,0);
\draw[thick,fill] (12,0) circle (0.3);
\draw[thick,bend left,distance=40] (12,0) to (15,0);
\draw[thick,fill] (15,0) circle (0.3);
\draw[thick,bend left,distance=40] (15,0) to (18,0);
\draw[thick,fill] (18,0) circle (0.3);
\draw[thick,bend left,distance=40] (18,0) to (21,0);
\draw[thick,fill] (21,0) circle (0.3);
\draw[thick,dotted,bend right,distance=25] (21,0) to
node[midway] {
}
(28,10);
\end{tikzpicture}
    \caption{$\widetilde\fT^{\nu}_\chi(\lambda,\mu)$}
    \label{fig:tensor}
\end{figure}

\begin{Remark}\label{rem:envelope}
  定理~\ref{thm:quiver} と定理~\ref{thm:tensor} を合わせると
  \begin{equation}\label{eq:92}
    H_{\mathrm{top}}(\widetilde{\fT}^{\nu}_\chi(\lambda,\mu))
    \cong\bigoplus_{\mu=\mu^1+\mu^2}
    H_{\mathrm{top}}(\fL_\chi(\lambda^1,\mu^1))\otimes
    H_{\mathrm{top}}(\fL_\chi(\lambda^2,\mu^2))
  \end{equation}
  という表現の同型が存在することが従う．しかし，上の既約成分の一対一対
  応が与える線形写像は，一般に $\mathfrak g$ の表現の同型を与えない．また，
  テンソル積表現は既約ではないので，表現の同型の取り方は一意ではない．
  有限型の箙に付随した箙多様体の場合は，一つのテンソル因子を最低ウェイト表現とみなすことにより，表現の同型写像を指定することができた\cite[Th.~5.9]{Na-Tensor}．
  Maulik-Okounkov \cite{2012arXiv1211.1287M} によって導入された stable
  envelope の理論は，一般の箙多様体において，標準的な表現の同型を幾何学的に与え
  る．
\end{Remark}

\begin{Remark}\label{rem:affineA}
\begin{NB}
  この構成は，$W$ をより多くの直和因子，たとえば $W = W^1\oplus W^2\oplus
  W^3$ と分解したときにも，ただちに一般化される．この例は，三つのテンソル積
  $V(\lambda^1)\otimes V(\lambda^2)\otimes V(\lambda^3)$に対応する．
\end{NB}%
  アファイン $A$ 型のディンキン図式(より一般にディンキン図式がループを
  持つ場合)に付随した箙多様体の場合は，別のタイプ
  の $\fM_0(\lambda,\mu)$, $\fM_\chi(\lambda,\mu)$ に作用する1パラメー
  タ部分群がある．頂点を $0$から $n-1$ に順番に $\ZZ/n\ZZ$ により番号を
  付ける．このとき $\Hom(V_{n-1},V_0)\oplus \Hom(V_{n-1},V_0)^*$ の部分
  をスカラー倍することが，誘導する $\fM_0(\lambda,\mu)$,
  $\fM_\chi(\lambda,\mu)$ への作用である．
  \begin{NB}
    When the diagram has no loop, this action on $\bN$ can be absorved
    by the $G$-action, hence is trivial on quotient spaces. But this
    is nontrivial when the diagram contains a loop.
  \end{NB}%
  この作用に関する $\fM_\chi(\lambda,\mu)$ の固定点は，自然
  に $A_\infty$ 型の箙多様体になる．頂点の番号づけ
  を $\ZZ/n\ZZ$ から $\ZZ$に変えることにより，$A_\infty$型の箙多様体に
  なるが，番号づけを $\ZZ$ から $\ZZ/n\ZZ$に変えることで与えられ
  る $A_\infty$ 型箙多様体から$A_{n-1}^{(1)}$ 型箙多様体への写像が，固
  定点集合の埋め込みに対応する．

  この作用に関して attracting setを考えると，リー環の準同
  型$\algsl_{n,\mathrm{aff}}\to \widehat{\mathfrak{gl}}(\infty)$ が誘導
  する表現が実現される．ここで後者は，無限サイズの行
  列 $(a_{ij})_{i,j\in\ZZ}$ ($a_{ij} = 0$ for $|i-j|\gg 0$) の全体の中
  心拡大のリー環である．
\end{Remark}

\begin{NB}
\subsection{層を用いた定式化}

定理~\ref{thm:quiver}と定理~\ref{thm:satake} をより近く比較するために，上の定義を偏屈層を用いて再定式化しよう．層を用いた合成積の定式化に馴染みのない読者は，この小節と次の小節は飛ばして構わない．

\cite[Cor.~10.11]{Na-alg} により(値域を $\pi$ の像に取り替えると)
$\pi\colon \fM_\chi(\lambda,\mu)\to \fM_0(\lambda,\mu)$ は semi-small
であることが知られている．従って，直像
$\pi_!(\mathcal C_{\fM_\chi(\lambda,\mu)})$ は半単純な偏屈層である．
ここで $\mathcal C_X$ は次元 $\dim X$ だけシフトした定数層，すなわち $\CC_X[\dim X]$ である．
点 $0$ は stratum をなし，
\begin{NB2}
\begin{equation*}
  \pi_!(\mathcal C_{\fM_\chi(\lambda,\mu)})\cong
  H_{\mathrm{top}}(\fL_\chi(\lambda,\mu)) \otimes \CC_{0}
  \oplus (\text{other summands}).
\end{equation*}
\end{NB2}%
\begin{equation*}
  H_{\mathrm{top}}(\fL_\chi(\lambda,\mu)) \cong
  \Hom(\CC_{0}, \pi_!(\mathcal C_{\fM_\chi(\lambda,\mu)})),
\end{equation*}
が成り立つ．ここで右辺の $\Hom$ は$\fM_0(\lambda,\mu)$上の偏屈層のなすアーベル圏における射の全体である．
(ある自然なstratificationに沿って局所定数層なものしか現れないので，そこに制限してもよい．)

上で注意したように $\pi$ は $\fM_{\chi_i}(\lambda,\mu)$ を経由する．
その写像を $\pi'$, $\pi''$ で表し，
$\pi = \pi'\circ\pi''$ と分解しよう．すると $\pi''$ も semi-small であり，
$\pi''_!(\mathcal C_{\fM_\chi(\lambda,\mu)})$ は，$\fM_{\chi_i}(\lambda,\mu)$ 上の半単純偏屈層である．
さらに，strata上の交叉コホモロジー複体の直和に同型である．
 (\cite[Prop.~15.3.2]{Na-qaff}の議論により，非自明な局所系に付随した
 交叉コホモロジーが生じないことが示される．) これを
\begin{equation*}
  \pi''_!(\mathcal C_{\fM_\chi(\lambda,\mu)})
  \cong \bigoplus_\alpha
  L_\alpha\otimes \mathrm{IC}(\fM_{\chi_i}^\alpha)
\end{equation*}
と表そう．ここで，$\fM_{\chi_i}^\alpha$ は strata，$L_\alpha$ は
$\fM_{\chi_i}^\alpha$ の点の上の $\pi"$ のファイバーの最高次のホモロジーである．
さらに stratum は，ある非負整数 $k$ を用いて
$\fM_{\chi_i}^{\mathrm{s}}(\lambda,\mu+k\alpha_i)$ と書けることが知られている．ここで，
添字 `s' は，安定な点のなす開集合を表している．\cite[Prop.~2.25]{Na-branching}を参照．

$\pi=\pi'\circ\pi''$であるから
\begin{equation*}
  H_{\mathrm{top}}(\fL_\chi(\lambda,\mu)) \cong
  \bigoplus_\alpha L_\alpha\otimes
  \Hom(\CC_0, \pi'_!\mathrm{IC}(\fM_{\chi_i}^\alpha))
\end{equation*}
となる．stratification $\fM_{\chi_i}(\lambda,\mu)=\bigsqcup \fM_{\chi_i}^\alpha$
は，閉埋め込み $\fM_{\chi_i}(\lambda,\mu)\hookrightarrow \fM_{\chi_i}(\lambda,\mu-\alpha_i)$ と整合的なので，
$\fM_{\chi_i}(\lambda,\mu-\alpha_i)$に新しく現れる stratum に対しては
 $L_\alpha = 0$ と解釈することにより，同じ $\alpha$ の集合の和として\begin{equation*}
  H_{\mathrm{top}}(\fL_\chi(\lambda,\mu-\alpha_i)) \cong
  \bigoplus_\alpha L_\alpha'\otimes
  \Hom(\CC_0, \pi'_!\mathrm{IC}(\fM_{\chi_i}^\alpha))
\end{equation*}
が成り立っているものとしてよい．

上の $e_i$, $f_i$ の定義の仕方により，分解を保つことが従う．つまり，
$e_i$, $f_i$は
\(
   L_\alpha \leftrightarrows L_\alpha'
\)
に定義される作用素と，
\(
   \Hom(\CC_0, \pi'_!\mathrm{IC}(\fM_{\chi_i}^\alpha)).
\)
の恒等写像のテンソル積である．

箙多様体の局所的な記述 (\cite[\S3]{Na-qaff}を参照)により，
 stratum $\fM_{\chi_i}^\alpha$ の点の上のファイバーは，ある別の箙多様体の
 $0$ のファイバー$\fL_\chi$ と同型であることが一般的に知られているが，
 特に今の状況では，$A_1$型の箙多様体になる．
$e_i$, $f_i$ の定義はこの記述と整合的である．すなわち，
\(
   L_\alpha \leftrightarrows L_\alpha'
\)
は，構成を $A_1$ 型の箙多様体の場合に行なったものを，上の記述を用いて
一般の箙多様体に移したものに他ならない．
$A_1$ 型の箙多様体は，グラスマン多様体の余接束と同型であり \cite[\S7]{Na-quiver}，その場合は
構成は，Ginzburg \cite{MR1111326} がその以前に与えたものに他ならない．

上の説明から明らかに，
$\Hom(\CC_0, \pi'_!\mathrm{IC}(\fM_{\chi_i}^\alpha))$ の次元は，
制限
\begin{equation*}
  V(\lambda)\!
  \downarrow^{\mathfrak g}_{\mathfrak{l}_i}
\end{equation*}
の既約表現の重複度に等しい．ここで， $\mathfrak l_i$ は $i$ に付随したLevi部分代数である．
\cite{Na-branching}を参照．
今の状況では $\mathfrak{l}_i$ は$\algsl_2$ とabelianなリー環の直和である．
添字 $\alpha$ は， $\mathfrak{l}_i$ の既約表現に対応する．$\fM_{\chi_i}^{\mathrm{s}}(V^0,W)$ が空でない
ことから
$\dim W_i + \dim V_{i-1}+\dim V_{i+1} \ge 2 \dim V_i$ が成り立ち，従って $\mathfrak l_i$ に対して支配的であることに注意しよう．

\subsection{テンソル積の層を用いた定式化}\label{subsec:sheaf_tensor}

層を用いた定式化を続けよう．次は \ref{thm:tensor}の再定式化を行う．
この結果は \cite{tensor2} による．

図式
\begin{equation*}
  \fM_0(\lambda,\mu)^{\nu}
  \xleftarrow{p}
  \fT^{\nu}_0(\lambda,\mu)
  \xrightarrow{j} \fM_0(\lambda,\mu),
\end{equation*}
を考える．ただし $j$ は埋め込みであり， $p$ は極限を与える
写像 $\lim_{t\to 0}\nu(t)\cdot$ である．双曲制限関手
 $\Phi$ を $p_* j^!$ により定義する．その基本的な性質は \cite{Braden,MR3200429} において議論されている．
 特に $\Phi$ は
$\pi_!(\mathcal C_{\fM_\chi(\lambda,\mu)})$ を単純偏屈層のシフトの
直和に，構成可能層の導来圏の対象として同型であることが知られている．

さらに今の設定では
$\Phi\pi_!(\mathcal C_{\fM_\chi(\lambda,\mu)})$ は半単純偏屈層
であること，すなわちシフトは必要ないことが知られている．
これは $\Phi$ が\cite{2014arXiv1406.2381B} の意味
で\emph{双曲 semi-small}\footnote{双曲semi-smallの定義
  は，\cite{2014arXiv1406.2381B}で与えられたが，実質的には\cite{MV2}に
  すでに現れている．固有射がsemi-smallであると押し出し写像は半単純偏屈
  層を半単純偏屈層に写すが，双曲制限関手において対応する概念である．}で
あることの帰結である．一方，双曲 semi-small であることは定
理~\ref{thm:tensor}(1) の帰結である．
最高次ホモロジーと双曲制限は
\begin{equation*}
  H_{\mathrm{top}}(\widetilde{\fT}_\chi^{\nu}(\lambda,\mu))
  \cong \Hom(\CC_0,\Phi\pi_!(\mathcal C_{\fM_\chi(\lambda,\mu)}))
\end{equation*}
により関係している．\cite[the paragraph after Lemma~4]{tensor2}を参照せよ．

定理~\ref{thm:tensor}  の主張(2) は，双曲制限
$\Phi\pi_!(\mathcal C_{\fM_\chi(\lambda,\mu)})$ の計算と理解することができる．
和で与える射
\begin{equation*}
  \sigma\colon \bigsqcup_{\mu=\mu^1+\mu^2}
  \fM_0(\lambda^1,\mu^1)\times\fM_0(\lambda^2,\mu^2)
  \to \fM_0(\lambda,\mu)^{\nu}
\end{equation*}
は，有限で全射であることがあることが知られている．
\cite[Lemma~1]{tensor2}.

\begin{Lemma}[\protect{\cite[Lemma~3]{tensor2}}]
  同型
  \begin{equation*}
    \sigma_! \bigoplus_{\mu=\mu^1+\mu^2}
    (\pi\times\pi)_! \mathcal C_{\fM_\chi(\lambda^1,\mu^1)\times
      \fM_\chi(\lambda^2,\mu^2)} \cong 
    \Phi\pi_!(
    \mathcal C_{\fM_\chi(\lambda,\mu)})
  \end{equation*}
  が存在する．
\end{Lemma}

これは，\subsecref{subsec:tensor1}における
$\widetilde{\fT}_\chi^{\nu}(\lambda,\mu)$ の記述からの帰結である．

$0$における茎をとることにより，同型 \eqref{eq:92}が従う．すで
に \ref{rem:envelope} で言及したように， \eqref{eq:92} は一意的ではない．
これは，上の補題の同型でも同様である．stable envelope を用いると，標準
的な同型を与えることができる．これは \cite[Lemma~4]{tensor2} の帰結であ
る．
\end{NB}%

\section{Kac-Moodyリー環に対する幾何学的佐武対応}\label{sec:Coulomb_Kac_Moody}

\subsection{クーロン枝と可積分最高ウェイト表現}\label{subsec:integrable2}

\secref{sec:quiver}の箙ゲージ理論のクーロン枝 $\cM(\lambda,\mu)$ に戻ろう．\subsecref{subsec:quiver_integrable}と同様に，頂点を自分自身に結ぶループを持たないと仮定する．
対応するKac-Moodyリー環を$\mathfrak g$とする．
注~\ref{rem:symmetrizable}で説明したように，クーロン枝の定義を修正することにより，
$\mathfrak g$が対称とは限らない場合にも拡張しておく．

\subsecref{subsec:quantum}で説明したように，クーロン枝は量子化から
誘導されるポアソン括弧と，非特異部分の上にシンプレクティック形式を持つ．
一般に，クーロン枝の特異点はBeauvilleの意味でのsymplecticな特異点である
と期待されている．これは箙ゲージ理論のクーロン枝の場合には正し
い\cite{2020arXiv200501702W}．特に，$\cM(\lambda,\mu)$は有限個
のsymplectic leaf からなる自然なstratification を持つ．

箙$Q$が有限型の場合には，このstratificationは
$\cM(\kappa,\mu)$ の非特異集合 $\cM^{\mathrm{s}}(\kappa,\mu)$ を用いて
\begin{equation*}
   \cM(\lambda,\mu) = \bigsqcup \cM^{\mathrm{s}}(\kappa,\mu)
\end{equation*}
で与えられる．ここで $\kappa$ は $\lambda\ge\kappa\ge\mu$ を満たす支配的余ウェイトを走る\cite{2019arXiv190209771M}．
$Q$がアファイン型のときには次のようになることが予想されており，
アファイン$A$型のときには正しいことが示されている\cite{2016arXiv160602002N}．
\begin{equation}\label{eq:94}
  \cM(\lambda,\mu) = \bigsqcup_{\kappa,\underline{k}}
  \cM^{\mathrm{s}}(\kappa,\mu)\times S^{\underline{k}}(\CC^2\setminus \{0\}/
  (\ZZ/\ell\ZZ))
\end{equation}
ここで， $\underline{k} = [k_1,k_2,\dots]$ は分割であ
り，
$\kappa$は$\lambda-|\underline{k}|c\ge\kappa\ge\mu$を満たす余ウェイトで
ある．(ここで，$c$ はアファイン・リー環の標準的な中心元である．) 簡単の
ために$\lambda$, $\mu$ のレベル$\ell$が $1$ より大きいと仮定しているが，レベ
ルが$1$のときの記述もほぼ同様である．

一般の箙ゲージ理論のクーロン枝のstratificationの場合に，\eqref{eq:94}の
対称積の一般化として，どんなものが生じるかは分かっていないが，ヒッグス枝
$\fM_0(\lambda,\mu)$のstratificationは\cite{CB}によって決定されており，
strataは一対一に対応しているのではないだろうかと，期待されている．

\subsecref{subsec:quiver_integrable}において，行列式の積で与えられる指
標 $\chi\colon\bG\to \CC^\times$ を選んでいたことを思い出そう．誘導する準同型
$\pi_1(\chi)\colon \pi_1(\bG)\to \pi_1(\CC^\times)$ と，
そのポントリャーギン双対
$\pi_1(\chi)^\wedge \colon \pi_1(\CC^\times)^\wedge \to
\pi_1(\bG)^\wedge$を考える．
\subsecref{subsec:hamtori} により，$\pi_1(\bG)^\wedge$ は $\cM(\lambda,\mu)$に作用する．
箙ゲージ理論の場合には $\pi_1(\bG)^\wedge = (\CC^\times)^n$である．(簡単のため，
すべての $i$ について$V_i\neq 0$と仮定する．)
同型 $\pi_1(\CC^\times)^\wedge = \CC^\times$を通じて，
$\pi_1(\chi)^\wedge$ を
$\pi_1(\bG)^\wedge = (\CC^\times)^n$内の1パラメータ変換群とみなし，
さらに$\pi_1(\chi)^\wedge$を$\chi$ で簡潔に表すことにする．

次の結果は，有限型のときは \cite{2017arXiv170900391K}，アファイン$A$型
のときは \cite[Prop.~7.30]{2016arXiv160602002N}
(\cite[Prop.~4.1]{2018arXiv181004293N}も参照)において証明され，一般の場
合にも成り立つと予想されている．
\begin{Conjecture}\label{lem:fixed}
  固定点集合$\cM(\lambda,\mu)^{\chi}$は，空集合か，もしくは一点である．
\end{Conjecture}

\eqref{eq:91}の類似として，次のattracting setを考えよう．
\begin{equation*}
  \fA_\chi(\lambda,\mu) \defeq
  \left\{ x\in \cM(\lambda,\mu) \,\middle|\,
    \text{$\lim_{t\to 0} \chi(t) x$が存在する}\right\}.
\end{equation*}
上の予想により，チルダ付きのバージョンと，なしのものは同じであることに注意しよう．

\begin{NB}
\begin{Remark}
  ヒッグス枝の場合，$\chi$ は特異点解
  消$\fM_\chi(\lambda,\mu)\to\fM_0(\lambda,\mu)$を与えた．同じ $\chi$
  がクーロン枝の場合に$\cM(\lambda,\mu)$に作用する1パラメータ変換群を
  与える．次の小節でみるように，$\fM_0(\lambda,\mu)$に作用する1パラメー
  タ変換群 $\nu$
  は，これとは反対に$\cM(\lambda,\mu)$の部分特異点解消(と変形)を与える
  ことになる．

  物理の文脈では，FIパラメータと質量パラメータの役割が，クーロン枝とヒッグス
  枝で入れ替わると説明されていたが，我々のクーロン枝の数学的な定義では
  この性質が数学的に厳密に確立されたわけである．
\end{Remark}
\end{NB}%

\cite{2016arXiv160403625B}で提唱されたKac-Moodyリー環に関す
る幾何学的佐武対応は，次の予想である．
\begin{Conjecture}\label{thm:satake}
  \textup{(1)} $\fA_\chi(\lambda,\mu)$ と各symplectic leafの共通部分は，
  空集合であるか，もしくはラグランジアン部分多様体である．

  \textup{(2)} 最高次のホモロジー群の直和
  \begin{equation*}
    \bigoplus_\mu H_{\operatorname{top}}(\fA_\chi(\lambda,\mu))
  \end{equation*}
  は，Kac-Moodyリー環 $\mathfrak
  g$ のLanglands双対$\mathfrak g^\vee$
  の，最高ウェイトが$\lambda$の可積分最高ウェイト表現の構造を持つ．
\end{Conjecture}

Kac-Moodyリー環 $\mathfrak g$ のLanglands双対 $\mathfrak
g^\vee$ は，$\mathfrak g$のカルタン行列を転置行列で取り替えて定まるもの
として定義する．

この予想は，$\mathfrak g$ が有限次元の場合は，通常の幾何学的佐武対応に
帰着させることにより，\cite{2017arXiv170900391K}で証明された．
ポイントは定理~\ref{thm:slice}のあとで注意したよう
に，$\fA_\chi(\lambda,\mu)$が，\eqref{eq:h:1}と同型であることである．
アファイン$A$型のときは\cite{2018arXiv181004293N}において証明された．

すでに言及したとおり，予想~\ref{thm:satake}と定理~\ref{thm:quiver}の主
張は，形式的にはよく似ている．この類似は生成元 $e_i$, $f_i$, $h$ の定義
の仕方にも見られる．カルタン部分環の元$h$ は次元ベクトルにより定められ
る．すなわち $H_{\mathrm{top}}(\fA_\chi(\lambda,\mu))$ が，ウェイト空間
$V_\mu(\lambda)$ に対応するように定められる．
一方$e_i$, $f_i$ は，定理~\ref{thm:quiver}の構成に現れた指標 $\chi_i$
に付随した多様体を調べることにより次のように定義される．$\chi_i$による
固定点集合
\(
  \cM(\lambda,\mu)^{\chi_i}
\)
と，その attracting set
\(
\fA_{\chi_i}(\lambda,\mu) \defeq \left\{ x\in\cM(\lambda,\mu)
  \,\middle|\, \text{$\lim_{t\to 0}\chi_i(t)x$が存在する}\right\}
\)
を考える．極限$\lim_{t\to 0}\chi(t)\cdot$を取ることにより，写像
$\fA_{\chi_i}(\lambda,\mu)\to \cM(\lambda,\mu)^{\chi_i}$が定まることに注意する．
さらに$\fA_\chi(\lambda,\mu)$内の$\chi_i$固定点集合
\(
   \fA_{\chi}(\lambda,\mu)^{\chi_i}
\)
を考える．これは，$\cM(\lambda,\mu)^{\chi_i}$ 内の$\chi$に関す
る attracting set とも同じである．
すると，もともとのattracting setは，ファイバー積
\begin{equation}\label{eq:93}
  \fA_\chi(\lambda,\mu) = \fA_\chi(\lambda,\mu)^{\chi_i}
  \times_{\cM(\lambda,\mu)^{\chi_i}} \fA_{\chi_i}(\lambda,\mu)
\end{equation}
として実現されることが分かる．

次に，$\cM(\lambda,\mu)^{\chi_i}$ は空集合か，もしくは$A_1$型の箙ゲージ
理論のクーロン枝 $\cM_{A_1}(\lambda',\mu')$ と同型であることが予想でき
る．ここで，$\mu' = \langle\mu,\alpha_i\rangle$ であり，$\lambda'$ は説
明しないが，ある余ウェイトである．したがっ
て，$\fA_\chi(\lambda,\mu)^{\chi_i}$ は$A_1$型の箙ゲージ理論のクーロン
枝のattracting setになる．

この記述が得られると，$e_i$, $f_i$ を $A_1$ すなわち $\algsl(2)$ の場合
に帰着させて定義することができる．この方法は，定理~\ref{thm:quiver}に
おける構成と形式的に似ている．
$e_i$, $f_i$ の正確な定義のためには，元々の幾何学的佐武対応のときと同じ
ように，$\cM(\lambda,\mu)$の偏屈層上の双曲
制限関手を使う必要があるが，この論説ではこれ以上の詳細には入らないものとする．

$\bigoplus_\mu
H_{\operatorname{top}}(\fA_\chi(\lambda,\mu))$は
$\fA_\chi(\lambda,\mu)$の既約成分の基本類が与える基底を持つことに注意し
よう．
上と同様の $A_1$ 型箙ゲージ理論への帰着により，既約成分の集
合$\bigsqcup \operatorname{Irr} \fA_\chi(\lambda,\mu)$に柏原のクリスタ
ルの構造が入り，さらにそれはKac-Moodyリー環 $\mathfrak
g^\vee$
の量子展開環$\mathbf U_q(\mathfrak g^\vee)$の可積分最高ウェイト表現の結
晶基底のクリスタルに同型になることが期待される．有限型のとき
は
\cite{2017arXiv170900391K}，アファイン$A$型のとき
は\cite[\S5(vi)]{2018arXiv181004293N}において証明されている．

\subsection{クーロン枝によるテンソル積}

\subsecref{subsec:tensor1}において，分解$W = W^1\oplus W^2$に付随し
た1パラメータ変換群$\nu\colon\CC^\times\to\GL(W)$を取ったこと
を思い起こそう．これは\subsecref{subsec:flavor}で説明したゲージ理論のフ
レーバー対称性の例と理解することができる．実際，上の分解と整合的
な$W_i$の基底を取り，$\GL(W_i)$の対角行列の全体 $T(W_i)$ を取っ
て，$\tilde\bG = \bG\times\prod_i T(W_i)$とする
と，\eqref{eq:7}の$\bN$は
$\tilde\bG$の表現に拡張される．そこで，
$\nu$ を $(\tilde\bG/\bG)^\vee$の指標と思って，
\subsecref{subsec:flavor}の構成を適用することができる．
従って，運動量写像のレベルを$\nu$にとる変形や，部分特異点解消を考え
ることができる．
ここでは，変形の方を考えて
$\cM^{\nu}(\lambda,\mu)$ で表そう．

\subsecref{subsec:tensor1}と同様に，$(\dim W^1)$,
$(\dim W^2)$に対応する支配的余ウェイトを $\lambda^1$, $\lambda^2$
で表そう．
予想~\ref{lem:fixed}と同様に次が期待される．
\begin{Conjecture}
  固定点集合 $\cM^{\nu}(\lambda,\mu)^\chi$ は，有限集合であり
  \(
  \bigsqcup_{\mu=\mu^1+\mu^2} \cM(\lambda^1,\mu^1)^\chi
  \times \cM(\lambda^2,\mu^2)^\chi
  \)
  との間に自然な一対一対応を持つ．
\end{Conjecture}

予想~\ref{lem:fixed}により $\cM(\lambda,\mu)^\chi$ は空集合か一点であり，
予想~\ref{thm:satake}と組み合わせれば，一点である必要十分条件
は$V_\mu(\lambda)\neq 0$である．従って，上の予想
は $\cM^{\nu}(\lambda,\mu)^\chi$ が空でない必要十分条件
は，$\mu$
がテンソル積$V(\lambda^1)\otimes V(\lambda^2)$のウェイトであることを主
張している．次のattracting setを考えよう．
\begin{equation*}
  \fA_\chi^{\nu}(\lambda,\mu) \defeq
  \left\{ x\in \cM^{\nu}(\lambda,\mu) \,\middle|\,
    \text{$\lim_{t\to 0} \chi(t) x$が存在する}\right\}.
\end{equation*}

次の予想は，表現のテンソル積の実現を与える．
\begin{Conjecture}
  $\fA_\chi^{\nu}(\lambda,\mu)$の最高次ホモロジー群の
  $\mu$に関する直和
  \begin{equation*}
    \bigoplus_\mu H_{\mathrm{top}}(\fA_\chi^{\nu}(\lambda,\mu))
  \end{equation*}
  は，$\mathfrak g^\vee$ の可積分表現の構造を持ち，それは
  テンソル積表現$V(\lambda^1)\otimes V(\lambda^2)$と同型である．
\end{Conjecture}

形式的には，この予想は定理~\ref{thm:tensor}の類似である．

予想~\ref{thm:satake}と同様に，有限型かアファイン$A$型のときは予想は正しい．
(\cite[Cor.~4.9]{2018arXiv181004293N}を参照せよ．)

アファイン$A$型の場合，注~\ref{rem:affineA}で説明した，ループが定め
る1パラメータ部分群をフレーバー対称性として使うことによ
り，$\cM(\lambda,\mu)$ の変形や部分特異点解消を考えることができる．上の
予想の類似として，変形においてattracting setを考えることによっ
て，$\algsl_{n,\mathrm{aff}}\to \widehat{\mathfrak{gl}}(\infty)$が誘導
する表現が実現される．これは技術的な注意であるが，生成元$e_i$, $f_i$が
定義関係式を満たすことの証明を，有限$A$型の場合に帰着させるために用いられる．

\begin{NB}
\subsection{層を用いた定式化}\label{subsec:sheaf}

\subsecref{subsec:sheaf_tensor}と同様に，固定点集合，attracting setと全
空間をつなぐ図式
\begin{equation*}
  \cM(\lambda,\mu)^{\chi}
  \xleftarrow{p}
  \fA_\chi(\lambda,\mu)
  \xrightarrow{j} \cM(\lambda,\mu),
\end{equation*}
を考え，双曲制限関手 $\Phi = p_* j^!$を定義しよう．もしも
$\cM(\lambda,\mu)^\chi = \emptyset$であったら，$\Phi = 0$と約束する．
このとき，$\Phi(\mathrm{IC}(\cM(\lambda,\mu))$は，$0$でなければ，
一点 $\cM(\lambda,\mu)^\chi$上の半単純偏屈層である．
これは\subsecref{subsec:sheaf_tensor}と同様に$\Phi$が双曲semi-smallであ
ることの帰結であり，一方双曲semi-smallであることは，予
想~\ref{thm:satake}(1)の帰結である．また，
\begin{equation*}
  H_{\mathrm{top}}(\fA_\chi(\lambda,\mu)) \cong
  \Phi(\mathrm{IC}(\cM(\lambda,\mu)))
\end{equation*}
となることも従う．

\subsecref{subsec:integrable2}の指標 $\chi_i$ を使って，$\fA_\chi(\lambda,\mu)$と$\fA_\chi(\lambda,\mu)^{\chi_i}$による双曲制限関手 $\Phi_i$ and $\Phi^i$ を考えよう．
ファイバー積の性質\eqref{eq:93}から
$\Phi = \Phi^i\circ \Phi_i$が従う．従って，
\( \Phi(\mathrm{IC}(\cM(\lambda,\mu))) =
\Phi^i(\Phi_i(\mathrm{IC}(\cM(\lambda,\mu))))  \)
が成り立つ．

$\cM(\lambda,\mu)^{\chi_i}$が$\cM_{A_1}(\lambda',\mu')$と同型であり，
symplectic leaves によるstratification が $\bigsqcup \cM_{A_1}^{\mathrm{s}}(\kappa',\mu')$ 
($\lambda'\ge\kappa'\ge \mu'$)であったことを思い出そう．
射 $p$, $j$は stratificationと整合的なので，双曲制限関手で現れる
偏屈層はstrataに沿って局所定数である．もう少し精密な議論を
使うと，非自明な局所系は現れないことも示すことができる．
従って，あるベクトル空間$M^{\lambda,\mu}_{\kappa',\mu'}$が存在して
\begin{equation*}
  \Phi_i(\mathrm{IC}(\cM(\lambda,\mu)))
  \cong \bigoplus_{\kappa'} M^{\lambda,\mu}_{\kappa',\mu'}
  \otimes\mathrm{IC}(\cM_{A_1}(\kappa',\mu'))
\end{equation*}
が成立する．

$\mu-\alpha_i$について同様にして，
\begin{equation*}
  \Phi_i(\mathrm{IC}(\cM(\lambda,\mu-\alpha_i)))
  \cong \bigoplus_{\kappa'} M^{\lambda,\mu-\alpha_i}_{\kappa',\mu'-2}
  \otimes\mathrm{IC}(\cM_{A_1}(\kappa',\mu'-2))
\end{equation*}
が成立する．ここで，$\cM(\lambda,\mu)$ と $\cM(\lambda,\mu-\alpha_i)$ の
古典的なクーロン枝の記述 (定理~\ref{thm:classical}) を用いると，
自然な同型$M^{\lambda,\mu}_{\kappa',\mu'}\cong
M^{\lambda,\mu-\alpha_i}_{\kappa',\mu'-2}$を構成することができる．
(\cite[Prop.~5.11]{2018arXiv181004293N}を参照せよ．) すると，
$\Phi(\mathrm{IC}(\cM(\lambda,\mu)))$と
$\Phi(\mathrm{IC}(\cM(\lambda,\mu-\alpha_i)))$との間に$e_i$, $f_i$を定義するためには
$A_1$型の場合に考察すれば十分であることが分かる．
$A_1$型のときは，クーロン枝やattracting setを具体的に記述することが可能で，
attracting set はベクトル空間になっていることが示される．この記述を用いて 
$e_i$, $f_i$ を定義することが可能である．
なお，同型$M^{\lambda,\mu}_{\kappa',\mu'}\cong
M^{\lambda,\mu-\alpha_i}_{\kappa',\mu'-2}$の構成は，有限次元の
通常の幾何学的佐武対応の場合には，$\bG_\cO$の作用を用いて与えられるが，
上の構成と同じ同型になっていることが分かる．従って，この構成は通常の
幾何学的佐武対応の一般化になっていることが従う．
\end{NB}%

{\bf 謝辞}\ 
本稿を執筆する機会を提供してくださった`数学'編集部の皆様に御礼を
申し上げます．
また，原稿を精読し多くの有意義なご意見と指摘をしていただいた
査読者の方々に感謝申し上げます．



\bibliographystyle{myamsalpha}
\bibliography{nakajima,mybib,coulomb,orthsymp}    


\end{document}